%% file: ppvArxiv.tex
\documentclass[twocolumn]{imsart}
\usepackage{imsart}
\usepackage{amsthm}
\usepackage{amsmath}
\usepackage{amsfonts}
\usepackage[numbers]{natbib}
\usepackage[colorlinks,citecolor=blue,urlcolor=blue,filecolor=blue,backref=page]{hyperref}
\usepackage{graphicx}
\usepackage{comment}
\usepackage{color}
\usepackage{subfigure}
\usepackage{fullpage}
\usepackage{widetext}
\usepackage{pstricks,pst-node,pst-text,pst-3d}

\startlocaldefs
\numberwithin{equation}{section}
\theoremstyle{plain}
\newtheorem{thm}{Theorem}[section]
\newtheorem{asmp}{Assumption}
\newtheorem{Lem}{Lemma}
\newtheorem{proposition}{Proposition}[section]
\newtheorem{cor}{Corollary}[section]
\newtheorem{exm}{Example}[section]
\newcommand{\upmodels}{\perp\!\!\!\perp}
\newcommand{\cind}[3]{#1 \upmodels #2\mid #3}
\newcommand{\mc}[1]{\mathcal{#1}}

\endlocaldefs

\begin{document}

\begin{frontmatter}
\title{Uncertainty Quantification Via the Posterior Predictive Variance}
\runtitle{On the Posterior Predictive Variance}

\begin{aug}

\author[addr1]{\fnms{Sanjay} \snm{Chaudhuri}\ead[label=e2]{schaudhuri2@unl.edu}},
\author[addr2]{\fnms{Dean} \snm{Dustin}\ead[label=e1]{deandust55@gmail.com}},
\author[addr3]{\fnms{Bertrand} \snm{Clarke}\ead[label=e3]{bclarke3@unl.edu}}

\runauthor{Chaudhuri et al.}
\address[addr1]{Sanjay Chaudhuri, Department of Statistics, University of Nebraska-Lincoln, NE, USA, 68583-0963
  \printead{e2}
}

\address[addr2]{Dean Dustin, First Citizens' Bank, Raleigh, NC, USA
    \printead{e1} 
}

\address[addr3]{Bertrand Clarke, Department of Statistics, University of Nebraska-Lincoln, NE, USA, 68583-0963
  \printead{e3}
}


\end{aug}

\begin{abstract}
We use the law of total variance to generate multiple expansions for the posterior predictive variance.   These expansions are sums of terms involving
conditional expectations and conditional variances and provide a quantification
of the sources of predictive uncertainty. Since 
the posterior predictive
variance is fixed given the model, it represents a constant quantity that is
conserved over these expansions.   The terms in the expansions can be
assessed in absolute or relative sense to understand the main contributors to the
length of prediction intervals.   We quantify the term-wise uncertainty
across expansions varying in the number of terms and the order of conditionates.
In particular, given that a specific term in one expansion is small or 
zero, we identify
the other terms in other expansions that must also be small or zero.
We illustrate this approach to predictive model assessment
in several well-known models.
\end{abstract}



\end{frontmatter}

\section{The Setting and Intuition}
\label{Intro}

Everyone uses prediction intervals (PI's) but few examine their structure
or more precisely how they should be interpreted in the context of 
a model with multiple components.  Often PI's seem overconfident (too narrow)
or useless (too wide).  More often than not, PI's are an afterthought to
modeling rather than the focus: most sample size selection procedures,
for instance, focus on estimation or testing, not prediction.

Both frequentist and Bayesian practitioners routinely report PI’s.
It is common for frequentists to estimate a model and then use it, perhaps 
even without adjustment, to give PI's; see \cite{Shen:etal:2004}
and more recently \cite{Bachoc:etal:2019}, \cite{Tian:2020}, and
\cite{Liang:etal:2025}.  By contrast,
it is common for Bayesians to give a PI simply by simulation from the posterior
predictive distribution and report the posterior predictive variance
(PPV) itself as a scalar summary without making
the relationship between the width of the PI and the components of the PPV explicit.  Here, 
instead, we treat the PPV as the quantity of interest.
The goal of this paper is therefore to pull together many existing
ideas about the PPV, and examples, in a relatively complete and organized way
so they can be seen as a coherent and informative body. 

Analogous to classical components-of-variance models, e.g., split plot designs,
we use the Law of Total Variance (LTV) to expand the PPV into interpretable contributions. 
We implicitly adopt a Bayesian
standpoint not because we accept it (although we do) but because of
its Containment Principle:  all relevant distributions exist and are
`contained' in 
a single measure space. 
Readers who prefer non-Bayesian formulations may regard our parameters simply as random variables.

Looking at a variance brings in the metric properties of 
a distribution.  The difference between using a PPV to get a PI and simulating from the posterior directly to get a
PI is that the former relies on posterior normality and the squared
error distances between random variables and their means -- essentially
regarding the mean as a constant random variable whereas the latter
only uses probabilities.  To see this, consider
a generic two-level example.  Suppose
\begin{eqnarray}
Z &\sim& w(z) \nonumber \\
Y^n = (Y_1, \ldots, Y_n ) &\sim& p(y \vert z)
\label{twolevelgen}
\end{eqnarray}
where $Y=Y^n$ is independent and identical distributed (IID)
data and $Z$ is a conditioning variable e.g., a unidimensional
parameter.  Denoting $n$ outcomes by $y^n = (y_1, \ldots , y_n)^T$,
the conditional
distribution $(Y_{n+1} \vert y^n)$ will give PI's.  Alternatively,
the PPV is ${\sf Var}(Y_{n+1} \vert y^n)$, and
the LTV gives
\begin{align}
{\sf Var}(Y_{n+1} \vert y^n) 
=& 
{\sf E} {\sf Var}_{Z\vert y^n} (Y_{n+1} \vert y^n, Z)\nonumber\\ 
+& 
{\sf Var}_{Z\vert y^n}{\sf E}(Y_{n+1} \vert y^n, Z)
\label{twotermgen}
\end{align}
which can be interpreted.
Loosely, the first term on the right in \eqref{twolevelgen}
is the variability from the 
likelihood and the second term on the right is the variability 
from $w$. Viewing predictive uncertainty via this expansion clarifies why
PI's have the width they do. 
In this case, if the conditional distributions $p(\cdot \vert z)$ 
concentrate tightly in the space of densities, 
the first term may be negligible, implying the hierarchy effectively has lower dimension than it appears. 
On the other hand, if the conditional means ${\sf E}(Y_{n+1} \vert y^n, Z)$ vary little across $Z$, the second term may be small.
We can inadvertently (or artificially) make a term small by 
choosing the densities in a parametric family to be close to each
other but not identical, essentially a metric property.

These ideas extend naturally to hierarchical models of arbitrary depth. 
Consider a hierarchical model (HM) for a response $Y=y$ given 
$\mc{Z} =  (Z_1, \ldots , Z_k, \ldots, Z_K)^T$ taking values 
$z = (z_1, \ldots , z_K)^T$
for some $K \in \mathbb{N}$:
\begin{eqnarray}
Z_1 &\sim& w(z_1) \nonumber \\
Z_2 &\sim& w(z_2 \vert z_1) \nonumber \\
\vdots & \vdots  &  \vdots  \nonumber \\
Z_K &\sim& w(z_K \vert z_1, \ldots , z_{K-1}) \nonumber \\
Y^n &\sim& p(y \vert z),
\label{hierarchicalmodelgeneric}
\end{eqnarray}
where the $w$'s represent distributions for the $Z_k$'s as indicated 
by their arguments and $p(\cdot \vert z)$ is the likelihood. 
Again, $n$ IID copies of $Y$ are denoted by $Y^n = (Y_1, \ldots , Y_n)^T$ 
with outcomes $y^n = (y_1, \ldots , y_n)^T$.   

Sequential application of the LTV to \eqref{hierarchicalmodelgeneric} 
produces $K+1$ variance components, each interpretable in terms of 
contributions from different levels of the hierarchy.  
Because \eqref{hierarchicalmodelgeneric} satisfies the Containment Principle, 
it is straightforward to see how assumptions on conditional distributions affect the 
expansion of the PPV.  
By tracking the relative sizes of these components, we obtain a common scale for all sources of uncertainty, allowing coherent comparison of contributions from different levels.  In particular, we will track how one term being zero in one 
sequence of uses of the LTV can imply how a related term is zero in 
another sequence of uses of the LTV.

This paper provides a conceptual framework for understanding 
predictive uncertainty via LTV expansions of the PPV.  
This framework applies broadly, offering insight into the components that 
determine PI length and the interpretation of predictive statements in
practice.
Our goal is not to introduce new mathematical theory or computational methods, 
but to organize and clarify existing ideas.  We focus on
showing how these expansions help interpret PI's and the structure 
of predictive inference.

\`A propos of this, it is well-recognized that prediction requires 
calibration and often
re-calibration; see \cite{Qian:etal:2025a} and \cite{Qian:etal:2025b}
and the references therein
for recent contributions in regression and classification respectively.
We dodge this question here because it is not immediately
germane to our analysis of the PPV.

This paper proceeds as follows.  In Sec. \ref{sec:2} we present
a series of examples that illustrate many of the properties of
expansions of the PPV for models like \eqref{hierarchicalmodelgeneric}.
In Sec. \ref{sec:quant}, we develop properties of the use of
LTV expansions for two and three-term cases.
In Sec. \ref{decomposition}, we discuss the use of
more general expansions for uncertainty quantification.
In Sec. \ref{Draper} we present some
computational work comparing how the terms in a three-term
expansion behaves as functions of its inputs
along with a data-driven example of this.
In Sec. \ref{discuss}, we discuss the implications and uses 
of these expansions for Uncertainty Quantification.  Details of derivations are
relegated to Appendices \ref{calcs3termnormal} and
\ref{calcs2wayANOVA}.

\section{Two-term and Three-term Expansions of Posterior Predictive Variance}
\label{sec:2}

Let $\mc{D}$ denote the available data, which includes $y^n$, and any other covariate that might be available. Given $\mc{D}$,  the posterior predictive density to future values $Y_{n+1}$, that is 
\begin{eqnarray}
Y \sim
p(y_{n+1} \vert \mc{D}) = \int p(y_{n+1} \vert v) w(v \vert \mc{D}) {\rm d} v,
\label{postpredfut}
\end{eqnarray}
where $w(v \vert \mc{D})$ is the posterior density. 
At this point, the PPV within the context of
the model \eqref{hierarchicalmodelgeneric}
is fixed.  Denote it by ${\sf Var}(Y_{n+1} \vert y^n)$.   When a random variable in the top
$K$ levels of the hierarchy are visible, we say it is explicit.  Otherwise, we say it is implicit.
Thus,  ${\sf Var}(Y_{n+1} \vert y^n)$ depends implicitly on the top $K$ levels of \eqref{hierarchicalmodelgeneric}.

\subsection{Two-term Expansions}
\label{sec:twoTerm}

Let $V_1\in\mc{Z}$.  From the Law of Total Variance (LTV) it immediately follows that: 
\begin{align}
  {\sf Var}(Y_{n+1} \vert {\cal{D}}) =& {\sf E}_{V_1 \vert {\cal{D}}} ( {\sf Var} ( Y_{n+1}\vert V_1, {\cal{D}}) )\nonumber\\
  +& {\sf Var}_{V_1 \vert {\cal{D}}} ({\sf E}(Y_{n+1} \vert V_1, {\cal{D}})).
\label{postpredvarV1}
\end{align}
We define the above expansion as a two-term expansion of the PPV conditional on $V_1$ and $\mc{D}$.  The expansion easily extends to subsets $V_1\subseteq\mc{Z}$ of size larger than one.  

Consider a few common examples of two-term expansions from the parametric Bayesian Hierarchical models.
\begin{exm}
\label{2termnorm}
Consider using the Law of Total Variance (LTV) 
on the posterior predictive variance (PPV) from a 
normal likelihood with a conjugate prior, i.e.,
\begin{align}
{\sf Var}(Y_{n+1} \vert y^n)
=&
{\sf E}_{\mu \vert y^n} {\sf Var}(Y_{n+1} \vert y^n, \mu)\nonumber\\
+&
{\sf Var}_{\mu \vert y^n}  {\sf E}(Y_{n+1} \vert y^n, \mu)
\label{LTV2norm}
\end{align}
where $\mu \sim N(\mu_0, \tau_0^2)$ and the
$Y_i$'s are independently and identically distributed 
(IID)
as ${\sf N}(\mu, \sigma_0^2)$, where $\mu_0$, $\sigma_0$ 
and $\tau_0$ are known.    Here, $y^n = (y_1, \ldots , y_n)^T$ is an outcome of the 
Vector $Y^n = (Y_1, \ldots, Y_n)^T$. 
It is easy to see that
$
p(y_{n+1} \vert y^n, \mu) = p(y_{n+1} \vert \mu),
$
where $\sigma_0$ has been suppressed in the notation.   
So, it is also easy to see that
$$
{\sf E}(Y_{n+1} \vert y^n, \mu) = \mu \quad \hbox{and} \quad 
{\sf Var}(Y_{n+1} \vert y^n, \mu) = \sigma_0^2.
$$
Since ${\sf Var}(Y_{n+1} \vert y^n, \mu)$ is a constant, its expectation under the posterior for $\mu$
is unchanged.  Thus, the first term on the right in \eqref{LTV2norm} is
$$
{\sf E}_{\mu \vert y^n} {\sf Var}(Y_{n+1} \vert \mu, y^n) = \sigma_0^2.
$$

For the second term on the right in \eqref{LTV2norm} recall
the posterior for $\mu$ given $y^n$ is 
${\sf N}(\mu_n, \tau_n^2)$
where 
$$
\mu_n = \left( \frac{n}{\sigma_0^2} + \frac{1}{\tau_0^2} \right)^{-1} \frac{n}{\sigma_0^2}\left(\bar{y} + \frac{\mu_0}{\tau_0^2} \right)
$$
and
$$
\tau_n^2 = \left( \frac{n}{\sigma_0^2} + \frac{1}{\tau_0^2} \right)^{-1} .
$$
Now,
\begin{eqnarray}
{\sf Var}_{\mu \vert y^n}( {\sf E}(Y_{n+1} \vert \mu, y^n) ) = {\sf Var}_{\mu \vert y^n} (\mu) =  \left( \frac{n}{\sigma_0^2} + \frac{1}{\tau_0^2} \right)^{-1} .
\nonumber
\end{eqnarray}
So, \eqref{LTV2norm} is
\begin{eqnarray}
{\sf Var}(Y_{n+1} \vert y^n) = \sigma_0^2 
+
\frac{\sigma_0^2\tau_0^2}{\sigma_0^2 + n \tau_0^2}
\nonumber
\label{normnorm}
\end{eqnarray}
in which the first term dominates.
The first term is the intrinsic variance of
$Y_{n+1}$ and the second term is the extra variability due to not knowing $\mu$.\hfill\qed
\end{exm}

\begin{exm}
\label{exm:expF}
It is easy to generalize \eqref{LTV2norm} to a one dimensional exponential family for $Y$, say 
$h(y) e^{\eta T(y) - \phi(\theta)}$ on some domain, equipped with its conjugate prior.  Indeed,
some routine calculations show
\begin{eqnarray}
{\sf Var}(Y_{n+1} \vert y^n)
=
{\sf E}_{\theta \vert y^n} \phi^{\prime \prime}(\theta) 
+
{\sf Var}_{\theta \vert y^n}  \phi^\prime(\theta).
\label{expfamgen}
\end{eqnarray}
Again, the first term represents intrinsic variability in $Y_{n+1}$
and the second term represents the extra uncertainty 
from not knowing $\theta$.\hfill\qed
\end{exm}

\begin{exm}
\label{exm:betabinom}
This applies, for instance, to the Beta-Binomial problem.
Let  $Y_i \sim Bin(m_i, p)$ where $\sum_{i=1}^n m_i = M$,
we set $S = \sum_{i=1}^n Y_i$,
and $p \sim {\sf Beta}(\alpha, \beta)$.  Conjugacy gives
that the posterior for $p$ is 
$(p \vert Y^n) \sim {\sf Beta} (\alpha_n, \beta_n)$
where $\alpha_n = \alpha+ S$ and $\beta_n = \beta+ M -S$.
Now, the posterior predictive distribution for
$(Y_{n+1} \vert Y^n)$ is, by direct calculation,
${\sf Beta-Binomial}(m_{n+1}, \alpha_n, \beta_n)$,
with variance
\begin{align}
{\sf Var}(Y_{n+1} \vert Y^n)&\nonumber\\ 
=& 
m_{n+1} 
\frac{\alpha_n \beta_n}{(\alpha_n + \beta_n)^2}
\frac{(\alpha_n + \beta_n + m_{n+1})}{(\alpha_n + \beta_n + 1)}.
\label{postvarBB}
\end{align}
We get the same result from evaluating the two terms in \eqref{expfamgen}.  For $\theta = \log(p/(1-p))$, the
canonical form is
$P(Y_{n+1} = y_{n+1} \vert \theta ) = C(m_{n+1}, y_{n+1})e^{y\theta - m_{n+1} \log (1 + \exp(\theta)}$,
so
$\phi(\theta) = m_{n+1} \log(1 + e^\theta)$.  
The general form of the conjugate prior is
$w(\theta \vert \alpha, \beta) \propto e^{\alpha \theta - \beta \phi(\theta)}$ giving the posterior
$w(\theta \vert y^n) \propto e^{(\alpha + s)\theta - (\beta + M) \log(1 +  \exp(\theta))}$.  Since 
$\phi^\prime(\theta) = m_{n+1} p$ and
$\phi^{\prime \prime}(\theta) = m_{n+1}p(1-p)$, we get
\begin{equation*}
{\sf E}_{\theta \vert y^n} \phi^{\prime \prime}(\theta) 
=
{\sf E}(m_{n+1}p(1-p) \vert y^n)\nonumber\\ 
=
m_{n+1} \frac{\alpha_n \beta_n}{\alpha_n + \beta_n +1}
\nonumber
\end{equation*}
and
\begin{align}
{\sf Var}_{\theta \vert y^n}  \phi^\prime(\theta)
=&
{\sf Var} (m_{n+1} p \vert y^n)\nonumber\\
=&
m_{n+1}^2 \frac{\alpha_n \beta_n}{(\alpha_n + \beta_n)^2(\alpha_n + \beta_n +1)}
\nonumber
\end{align}
which, upon summing and re-arranging, gives 
\eqref{postvarBB}.\hfill\qed
\end{exm}

\begin{exm}
In canonical exponential family form, 
the {\sf Poisson}($\lambda$) is
\begin{eqnarray}
P(Y_i = y_i\vert \theta ) 
= (1/y!)e^{y\theta - \exp(\theta)},
\nonumber
\end{eqnarray}
so $\phi(\theta) = e^\theta = \lambda$.
The general form of the conjugate prior in $\theta$
is $w(\theta) \propto e^{\alpha\theta - \beta \exp{\theta}}$
giving the posterior 
$w(\theta \vert y^n) \propto e^{(\alpha +S)\theta - (\beta + n)\exp(\theta)}$.
By direct calculation, for $s = \sum y_i$,
\begin{eqnarray}
{\sf Var}(Y_{n+1} \vert y^n) 
=
\frac{\alpha + s}{\beta +n} \left(1 + \frac{1}{\beta+n} \right).
\label{postvarPG}
\end{eqnarray}
Since $\phi^\prime(\theta) = \phi^{\prime \prime}(\theta)  = e^\theta$, we get that the PPV is
the sum of
\begin{eqnarray}
{\sf E}(\phi^{\prime\prime}(\theta) \vert y^n)
=
\frac{\alpha + s}{\beta + n}
\quad \hbox{and} \quad
{\sf Var}(\phi^\prime(\theta) \vert y^n) = \frac{\alpha + s}{(\beta+n)^2} 
\nonumber
\end{eqnarray}
that gives \eqref{postvarPG}.
The Poisson$(\lambda)$ distribution with a Gamma prior
(giving a negative binomial posterior)
can also be worked out explicitly.\hfill\qed
\end{exm}

\subsection{Three-term Expansions}
\label{sec:threeTerm}

As the three foregoing examples
indicate, our focus here is to develop 
and study expansions of
the PPV.
As will be seen, our goal is to find ways to reduce the
number of terms in multi-term
expansions to eliminate
unnecessary conditioning variables.

To see that this is possible, let $\{V_1,V_2\}\subseteq\mc{Z}$, and  
extend \eqref{LTV2norm} by a second use of the LTV in the first term, 
\begin{align}
{\sf Var} ( Y_{n+1}\vert V_1, {\cal{D}}) =& 
{\sf E}_{V_2 \vert V_1,{\cal{D}}} ( {\sf Var} ( Y_{n+1}\vert V_1, V_2, {\cal{D}}) )\nonumber\\
+& {\sf Var}_{V_2 \vert V_1, {\cal{D}}} (E (Y_{n+1} \vert V_1, V_2, {\cal{D}})).
\label{LTVgeniter}
\end{align}

By substituting the above expression in \eqref{postpredvarV1}, one gets the \emph{three-term} expansion of the PPV as:
\begin{subequations}
\label{LTV2termgen}
\begin{align}
 ~& {\sf Var}( Y_{n+1} \vert {\cal{D}})\nonumber\\
  =& {\sf E}_{V_1 \vert {\cal{D}} } {\sf E}_{V_2 \vert V_1,{\cal{D}}} ( {\sf Var} ( Y_{n+1}\vert V_1, V_2, {\cal{D}}) ) 
\label{eq:1a} \\
+&
{\sf E}_{V_1 \vert {\cal{D}} } {\sf Var}_{V_2 \vert V_1, {\cal{D}}} ({\sf E} (Y_{n+1} \vert V_1, V_2, {\cal{D}}))
\label{eq:1b}\\
+&
{\sf Var}_{V_1 \vert {\cal{D}}} ({\sf E} (Y_{n+1} \vert V_1, {\cal{D}})).
\label{eq:1c}
\end{align}
\end{subequations}
Note that we could have used a two-term
expansion for the PPV simply by using only one
of the two conditioning variables.  The order of conditioning, however, matters. The three-term expansion defined in \eqref{LTV2termgen} depends on the sequence in which we condition on $V_1$ and $V_2$. That is, by conditioning on $V_2$ first and then on $V_1$, we can alternatively write \eqref{LTV2termgen} as:

\begin{subequations}
\label{LTV2termgen2}
\begin{align}
~&{\sf Var}( Y_{n+1} \vert {\cal{D}})\nonumber\\ 
=&{\sf E}_{V_2 \vert {\cal{D}} } {\sf E}_{V_1 \vert V_2,{\cal{D}}} ( {\sf Var} ( Y_{n+1}\vert V_2, V_1, {\cal{D}}) )\label{eq:2a} \\
+&
{\sf E}_{V_2 \vert {\cal{D}} } {\sf Var}_{V_1 \vert V_2, {\cal{D}}} ({\sf E} (Y_{n+1} \vert V_2, V_1, {\cal{D}})) \label{eq:2b}\\
&+
\hbox{Var}_{V_2 \vert {\cal{D}}} (E (Y_{n+1} \vert V_2, {\cal{D}})).\label{eq:2c}
\end{align}
\end{subequations}

Terms \eqref{eq:1a} and \eqref{eq:2a} are equal in general. However, neither \eqref{eq:1b} and \eqref{eq:2b} nor \eqref{eq:1c} and \eqref{eq:2c} need be the same.  Furthermore, the fact that any of these terms is zero does not reduce the three-term expansion to a valid two-term expansion as defined in \eqref{postpredvarV1}. 

We consider some illustrative examples of three-term expansions below:
\begin{exm}
\label{exm:NNG} 
Again, consider a normal likelihood, but this time with
a normal prior on the mean and a Gamma prior on the 
variance.  Let 
$Y_i \sim {\sf N}(\mu, \/\lambda^2)$ be IID for 
$i=1, \ldots, n$ and use the conjugate priors
$\mu \sim {\sf } {\sf N}(\mu_0, 1/(\kappa_0 \lambda^2) )$ 
with 
$\lambda^2 \sim {\sf Gamma}(\alpha_0, \beta_0)$.
This gives two three-term expansions depending on 
whether we condition on $\mu$ first
(and then $\lambda$) or $\lambda$ first (and then $\mu$).  

Conditioning on $\mu$ first, i.e., setting $V_1 = \lambda^2$ 
and $V_2 = \mu$ in \eqref{LTV2termgen} we get
\begin{subequations}
\begin{align}
 ~&{\sf Var}(Y_{n+1} \vert y^n)\nonumber\\
 =& {\sf E}_{\lambda^2\vert y^n} {\sf E}_{\mu \vert y^n, \lambda^2 } {\sf Var}(Y_{n+1} \vert y^n, \mu, \lambda^2)
\label{EEVar}\\
+&
{\sf E}_{\lambda^2\vert y^n} {\sf Var}_{\mu \vert y^n, \lambda^2}{ \sf E}(Y_{n+1} \vert y^n, \mu, \lambda^2)
\label{EVarE}\\
+&
{\sf Var}_{\lambda^2\vert y^n} {\sf E}_{\mu \vert y^n, \lambda^2} {\sf E}(Y_{n+1} \vert y^n, \mu, \lambda^2) ,
\label{VarEE}
\end{align}
\end{subequations}
in which it is easy to see that
\eqref{VarEE} is zero:

\begin{align}
~&{\sf Var}_{\lambda^2\vert y^n} {\sf E}_{\mu \vert y^n, \lambda^2} {\sf E}(Y_{n+1} \vert y^n, \mu, \lambda^2) 
= {\sf Var}_{\lambda^2\vert y^n} {\sf E}_{\mu \vert y^n, \lambda^2} (\mu)\nonumber\\
=&{\sf Var}_{\lambda^2\vert y^n} \left( \frac{\kappa_0\mu_0 + n \bar{y}}{\kappa_0 +n}\right) = 0.\nonumber
\end{align}

For \eqref{EEVar} and \eqref{EVarE}, we use the fact that, by conjugacy, there is an
$\alpha_n$ and $\beta_n$ so  that $\lambda^2 \vert y^n \sim {\sf Gamma}(\alpha_n, \beta_n)$.
This gives that
$$
{\sf E}_{\lambda^2 \vert y^n} \left( \frac{1}{\lambda^2}  \right) = \frac{\beta_n}{\alpha_n -1 }. 
$$
Now, dropping the conditioning on $y^n$ in the variance on the right of \eqref{EEVar} it is
\begin{align}
{\sf E}_{\lambda^2\vert y^n} {\sf E}_{\mu \vert y^n, \lambda^2 } {\sf Var}(Y_{n+1} \vert \mu, \lambda^2)
=&
{\sf E}_{\lambda^2\vert y^n} {\sf E}_{\mu \vert y^n, \lambda^2 } \left( \frac{1}{\lambda^2} \right)\nonumber\\
=& \frac{\beta_n}{\alpha_n -1 }.
\label{Fubiniterm}
\end{align}
Likewise, we can show that for $\kappa_n = n+\kappa_0$, \eqref{EVarE} is
\begin{align}
~&{\sf E}_{\lambda^2\vert y^n} {\sf Var}_{\mu \vert y^n, \lambda^2 } {\sf E}(Y_{n+1} \vert \mu, \lambda^2)\nonumber\\
=&
{\sf E}_{\lambda^2\vert y^n} {\sf Var}_{\mu \vert y^n, \lambda^2 } (\mu)\nonumber\\
=&
{\sf E}_{\lambda^2\vert y^n}\left(  \frac{1}{\lambda^2 \kappa_n} \right) =  \frac{\beta_n}{\kappa_n(\alpha_n -1)}.
\end{align}
Thus, we have that
\begin{eqnarray}
{\sf Var}(Y_{n+1} \vert y^n) = \left( \frac{\kappa_n + 1}{\kappa_n} \right) \frac{\beta_n}{\alpha_n -1 }.
\label{sumof3terms}
\end{eqnarray}

If we condition on $\lambda$ first and then $\mu$ we find that
\begin{subequations}
\begin{align}
~&{\sf Var}(Y_{n+1} \vert y^n)\nonumber\\
=& {\sf E}_{\mu\vert y^n} {\sf E}_{\lambda^2 \vert y^n, \mu } {\sf Var}(Y_{n+1} \vert y^n, \mu, \lambda^2)
\label{EEVar1}\\
+&
{\sf E}_{\mu\vert y^n} {\sf Var}_{\lambda^2 \vert y^n, \mu}{\sf E}(Y_{n+1} \vert y^n, \mu, \lambda^2)
\label{EVarE1}\\
+&
{\sf Var}_{\mu\vert y^n}  {\sf E}(Y_{n+1} \vert y^n, \mu ) .
\label{VarEE1}
\end{align}
\end{subequations}
Parallel to \eqref{VarEE}, it is easy to see that \eqref{EVarE1} is zero.
By Fubini, \eqref{EEVar1} is the same as \eqref{EEVar} as given by \eqref{Fubiniterm}.
Finally, since ${\sf Var}(Y_{n+1} \vert y^n)$ is constant independent of the condition,
we can solve for \eqref{VarEE1}.  If desired, we can calculate ${\sf Var}(Y_{n+1} \vert y^n)$ directly
and hence verify \eqref{sumof3terms}.  We give one version of this in Appendix \ref{calcs3termnormal}.

Comparing the two orders of conditioning, we see that in both 
the EEVar terms are the same.
In the first expansion, the VarEE term is zero, whereas in 
the second expansion, 
the EVarE term is zero.   Finally, in the first, the EVarE 
term has the $\kappa_n$ while
in the second, the VarEE term has the $\kappa_n$.  Thus, it is not
a priori clear which terms will dominate in three-term expansions.
\hfill\qed
\end{exm}

We also see that the interpretation of the two-term expansion (intrinsic 
variability plus extra parameter uncertainty) generalizes to
the three term case.  The first term continues to represent 
the intrinsic variability, but the second and third terms summarize the
extra variability due to parameter uncertainty -- the second term
for the first conditioning parameter and the third term for the
second conditioning parameter.

\begin{exm}
\label{exm:BPG}
As a second three-term example, suppose the $Y_i$'s are IID ${\sf Bin}(N_i, p)$ where $p$ is fixed and the $N_i$'s are drawn IID from a 
${\sf Poisson}(\lambda)$ and $\lambda \sim {\sf Gamma}(a, b)$.  
It is not hard to verify that
\begin{eqnarray}
{\sf Var}(Y_{n+1} \vert y^n)
=
p \frac{\sum_{i=1}^n y_i + a}{b+pn} \left(1 + \frac{p}{b + pn}\right)
\label{BPG}
\end{eqnarray}
and, as in the normal case, there are two 3-term expansions depending
on the order of conditioning on the $N_i$'s and $\lambda$.
The more natural ordering conditions on $N$ first:
\begin{subequations}
\begin{align}
  ~&{\sf Var}(Y_{n+1} \vert y^n)\nonumber\\
=& {\sf E}_{\lambda\vert y^n} {\sf E}_{N \vert y^n, \lambda } {\sf Var}(Y_{n+1} \vert y^n, \lambda, N)
\label{EEVar1B}\\
+&
{\sf E}_{\lambda\vert y^n} {\sf Var}_{N \vert y^n, \lambda}{\sf E}(Y_{n+1} \vert y^n, \lambda, N)
\label{EVarE1B}\\
+&
{\sf Var}_{\lambda \vert y^n}  {\sf E}(Y_{n+1} \vert y^n, \lambda ) .
\label{VarEE1B}
\end{align}
\end{subequations}
Let $s=\sum y_i$.  The corresponding terms are
\begin{eqnarray}
p(1-p) \frac{s+a}{b+pn} + p^2\frac{s+a}{b+pn} + p^2\frac{s+a}{(b+pn)^2}
\label{functionofs}
\nonumber
\end{eqnarray}
that clearly sum to \eqref{BPG}.  In order, the terms represent binomial,
Poisson, and Gamma variability, and all terms are functions of $s$
with the last term being of smaller order in $n$ than the first two.  
If the model is believed,
one can plot curves for the four terms as a function of $s$
to see which, if any, are small enough to be ignored.

If we condition in the reverse order, the $N$'s can be 
integrated out to give a
{\sf Poisson-Gamma} model, i.e., it reduces to a two-term expansion
because the
term parallel to \eqref{EVarE1B} is zero.
We have that the PPV is
\begin{eqnarray}
{\sf E}_{N \vert y^n} {\sf E}_{\lambda \vert y^n, N } 
{\sf Var}(Y_{n+1} \vert y^n, \lambda, N)
+
{\sf Var}_{N \vert y^n}  {\sf E}(Y_{n+1} \vert y^n, N ) 
\nonumber \\
=
\left( p(1-p) \frac{s+a}{b+pn} \right)
+
\left( p^2\frac{s+a}{b+pn} + p^2\frac{s+a}{(b+pn)^2} \right) .
\nonumber
\end{eqnarray}
\end{exm}
If we redo the calculations with a single fixed $N$ chosen at
the beginning according to a ${\sf Gamma}(a, b)$, we get the same
PPV and terms.   The reason is that the terms only depend on $s$ and this
updates the Poisson rate $p\lambda$ the same way in both cases.

\begin{exm}Here, we see that the Bayes model average (BMA) can be seen
as a two or three-term hierarchical model.
Let $j= 1, \ldots, J$ index a collection of models ${\cal{M}} = \{ M_1, \ldots , M_J\}$.
Assume each $M_j$ consists of a likelihood $p(y \vert \theta_j)$ and a prior 
$w(\theta_j, j) = w(\theta_j \vert j )w(j)$ where the across models prior $w(j)$ is discrete.
Writing $J$ for $j$ as a random variable, as well as for the number of models, will cause no confusion
because the context will indicate which is meant.  Now,
we can represent this as a two-level hierarchical model
\begin{eqnarray}
(J, \theta_J) &\sim& w(\theta_j,j) \nonumber \\
Y &\sim& p(y \vert \theta_j).
\label{BMA2term}
\end{eqnarray}
Now, the $L^2$ BMA predictor is
\begin{equation}
E(Y_{n+1} \vert y^n) = \sum_{j=1}^J E(Y_{n+1} \vert y^n, M_j)W(M_j \vert y^n),
\label{L2BMApred}
\end{equation}
where $W(\cdot\vert y^n)$ is the discrete posterior probability.

In \eqref{L2BMApred},  the two conditioning random variables, namely $J$ and $\theta_j$ are treated explicitly and implicitly, respectively.   In this case, it is not hard to see that one usage of the LTV recovers the usual formula for the PPV.
Indeed,  using the expression for posterior variance from p. 383 of \cite{Hoeting:etal:1999}, we find that
\eqref{L2BMApred} is
\begin{align}
~&{\sf Var}(Y_{n+1} \vert y^n)\nonumber\\ 
=&\sum_{j=1}^J {\sf Var}(Y_{n+1} \vert y^n, M_j) W(M_j \vert y^n)
\nonumber \\
&+ \sum_{j=1}^J W(M_j \vert y^n)\left\{{\sf E}(Y_{n+1} \vert y^n, M_j) 
- {\sf E}(Y_{n+1} \vert y^n)\right\}^2 \nonumber \\
=&  {\sf E}{\sf Var}(Y_{n+1} \vert M_J, y^n)+{\sf Var} ({\sf E}(Y_{n+1} \vert M_J, y^n)),
\label{LTVBMAsimple}
\end{align}
ie. \eqref{LTVBMAsimple} is the result of using the LTV and conditioning on $M_k$.
We have treated $\theta_j$ implicitly by integrating over it,
before conditioning on the $M_j$'s.   Reversing this, i.e., integrating over $j$ and using the LTV
with $\Theta_k$'s would have been mathematically well-defined but statistically inappropriate for BMA. 

When the first term on the right EVar in \eqref{LTVBMAsimple} is large,  we 
see most variability is in the predictive distributions from the high posterior probability models rather than from the variability across models.
The second term on the right being small means that it doesn't matter very 
much which model 
you use for prediction.  On the other hand, if VarE is large, model 
selection is important
but the smallness of the E Var term means all the high posterior probability
models are good.

Now write \eqref{BMA2term} as an equivalent three-level hierarchical model:
\begin{eqnarray}
J &\sim& w(j) \nonumber \\
\theta_j \vert J=j&\sim& w(\theta_j \vert j) \nonumber \\
Y &\sim& p(y \vert \theta_j).
\label{BMA2termequiv}
\end{eqnarray}
Now, using \eqref{LTV2termgen} with the natural choices of $V_1$ being the models and $V_2$ being the parameters we get:
\begin{align}
~&{\sf Var}(Y_{n+1} \vert y^n)\nonumber\\ 
=&  {\sf E}_J {\sf E}_{\Theta_J \vert y^n, M_J} {\sf Var}_{Y_{n+1} \vert y^n, M_J \Theta_j} (Y_{n+1} \vert  y^n, M_J,  \Theta_J) 
\nonumber \\
+&{\sf E}_J {\sf Var}_{\Theta_J \vert y^n, M_J} {\sf E}_{Y_{n+1} \vert  y^n, M_J, \Theta_J}(Y_{n+1} \vert y^n, M_J,  \Theta_J ) \nonumber \\
+&
{\sf Var}_J ({\sf E}(Y_{n+1} \vert M_J, y^n)).
\label{LTVBMA3term}
\end{align}
In this treatment of PPV, the relative size of the terms is a tradeoff among
the size of the model list,  the proximity of the parametric models on the list to each
other,  the across-models prior weights on models on the list, and the within-model priors.
\end{exm}

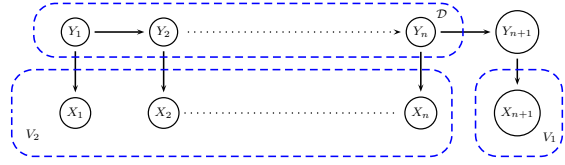
\begin{figure}[t]
  \resizebox{.8\columnwidth}{!}{\input{hmmPic.tex}}
  \smallskip
  \caption{Schematic diagram of a Hidden Markov model. With the definitions of $V_1$, $V_2$ and $\mc{D}$ as indicated, the conditional independence relationships $\cind{V_1}{V_2}{\mc{D}}$ and $\cind{Y_{n+1}}{V_2}{(V_2,\mc{D})}$ hold.}
  \label{fig:hmm}
\end{figure}

Although the examples of HM's we have seen so far are `vertical'
in the sense that each $Z_{k}$ in \eqref{hierarchicalmodelgeneric}
sits `above' $Z_{k+1}$, this is not necessary for our expansions.
Indeed, consider the following `horizontal' example diagrammed
in Fig. \ref{fig:hmm}.

\begin{exm}
\label{exm:hmm} 
Consider a Hidden Markov model where we observe
$X_1, \ldots, X_n$ assumed to be generated from
the hidden outcomes of a Markov process $Y_1, \dots , Y_n$
respectively.  The problem is
to predict $Y_{n+1}$ using the earlier $Y_i$'s as our `data'
${\cal{D}}$ even though they are unobserved. 
Since the $V_k$'s in \eqref{hierarchicalmodelgeneric} 
are simply random variables, not necessarily parameters,
we can use them in a three-term PPV expansion 
for ${\sf Var}(Y_{n+1} \vert y^n)$.

If we use $X_1, \ldots, X_n$ as $V_1$ 
and $X_{n+1}$ as $V_2$, then by construction, each $X_i$ depends 
only on its $Y_i$ for $i=1, \ldots, n+1$.  Now, suppose we have predicted values $\hat{Y}_i$ of $Y_i$
that are functions of the $X_1$, $X_2$, $\ldots$, $X_i$ and we use $\widehat{\mc{D}}=\left\{\hat{Y}_1,\hat{Y}_2,\ldots,\hat{Y}_n\right\}$.  
Our goal here is to find 
$\widehat{{\sf Var}}(Y_{n+1}\vert\widehat{\mc{D}} )
= {\sf Var} [Y_{n+1}\vert \hat{Y}_1, \hat{Y}_2, \ldots, \hat{Y}_{n}]$. 
This is valid from a prediction standpoint. 

With our definition of $V_1$ and $V_2$ a three-term expansion of $\widehat{{\sf Var}}(Y_{n+1}\vert \widehat{\mc{D}} )$ readily follows. We will show later that 
this expansion reduces to a two-term expansion.\qed
\end{exm}

Here, we limit ourselves to iteratively using the LTV on the EVar term,
each iteration bringing in one more $Z_k$.
There is nothing to stop us from applying the LTV in the VarE term as well;
however, such terms are very difficult to handle.
On the other hand, the PPV is constant over all these expansions and so
represents a `conservation of variance law'.  
Thus, while the full `LTV-scope' of a PPV contains 
many terms from using the LTV in all possible ways, we focus
on the subset of the LTV-scope where each term has exactly one 
variance operation
that moves from left to right with appropriate conditioning
and arises from using the LTV in the first, ie. EVar, term.   
We call this the Cochran Scope
or C-Scope for short.  This terminology recognizes the analogy
between a set of LTV expansions and the sums of squares decompositions 
that arise in frequentist ANOVA; see \cite{Dustin:Ghosh:Clarke:2025}.
Henceforth, we limit our attention to the C-scope's of HM's.
Proposition \ref{Cscopecard1} counts the number of expansions as a function
of $K$.

\section{Uncertainty Quantification for Two and Three-term Expansions}
\label{sec:quant}

In this section, we look at the behaviour of the individual terms 
in two- and three-term expansions.  Our goal is to identify
conditional independence assumptions -- that we call structural --
under which some of the terms in the RHS of equations 
\eqref{postpredvarV1} and \eqref{LTV2termgen} are small, perhaps zero.  
One implicit goal is to see if a level of the HM can be dropped,
as a consequence of dropping terms in the expansion.
For simplicity, we assume that both $V_1$ and $V_2$ are univariate.
It will be seen that these results generalize.  Indeed, three-term expansions 
exhibit all possible generic characteristics of expansions
for $K$ level HM's.


By an application of Fubini's Theorem, it is easy to see that
terms \eqref{eq:1a} and \eqref{eq:2a} are equal, i.e,
for any $V_1 ~ \mbox{and} ~ V_2 \in \mc{Z}$
\begin{align}
  ~&{\sf E}_{V_1 \vert {\cal{D}} } {\sf E}_{V_2 \vert V_1,{\cal{D}}} ( {\sf Var} ( Y_{n+1}\vert V_1, V_2, {\cal{D}}) )\nonumber\\
  =&{\sf E}_{V_2 \vert {\cal{D}} } {\sf E}_{V_1 \vert V_2,{\cal{D}}} ( {\sf Var} ( Y_{n+1}\vert V_2, V_1, {\cal{D}}) )\nonumber\\
    =&
    {\sf E}_{V_1,V_2 \vert {\cal{D}} } ( {\sf Var} ( Y_{n+1}\vert V_2, V_1, {\cal{D}}) ).
\end{align}
Indeed, in any valid modeling setting,  this term will be strictly
positive and leading in the sense that, typically,
no other term will be asymptotically larger
as $n$ increases.

\subsection{Uncertainty Quantification under Structural Assumptions}
\label{condlindpdce}

Here we give a sufficient condition for
the reduction of three-term expansions two term
expansions.  Then, we also show that, at least in the normal case, 
the higher the level in
the hierarchical model, the larger the conditional variances are.

\begin{asmp}
\label{asmp:condInd}
$Y_{n+1}$ and $\cal{D}$ are conditionally 
independent of $V_2$ given $V_1$, i.e.,
\begin{equation}
\label{condInd}
\cind{\left(Y_{n+1},\cal{D}\right)}{V_2}{V_1}.
\end{equation}
\end{asmp}

\begin{figure}[t]
\resizebox{.8\columnwidth}{!}{\input{bhm.tex}}
  \smallskip
\caption{Diagram for a `vertical' three-level HM.  It is easy to think
of $V_1$ is a parameter and $V_2$
as a hyperparameter.  Further conditioning, e.g., by including a $V_3$ 
would extend the diagram upwards. }
  \label{fig:horizontal}
\end{figure}
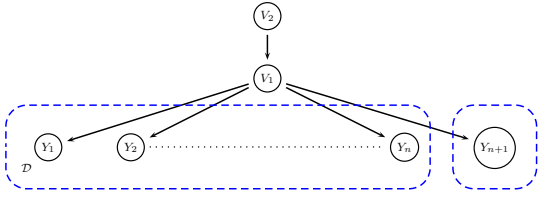

We call this structural because it must hold for all values of $V_1$ and $V_2$.  
That is, the condition \eqref{condInd} will be satisfied in any 
HM with $V_1$ being a parameter and $V_2$ being a hyperparameter; see Fig. \ref{fig:horizontal}.

\begin{thm}
\label{thm:struct}

Suppose Assumption \ref{asmp:condInd} holds.  Then, the term \eqref{eq:1b} 
in the three-term expansion \eqref{LTV2termgen} is zero.  
Furthermore, the expansion \eqref{LTV2termgen} reduces to 
a two-term expansion conditional on $V_1$ and $\cal{D}$.

\begin{proof} 
Assumption \ref{asmp:condInd} is equivalent to (see \citep{lauritzenBook})
$\cind{Y_{n+1}}{V_2}{\left(V_1,\cal{D}\right)}$ and 
$\cind{\mathcal{D}}{V_2}{V_1}$ taken together.  
From the first we get
    \begin{equation*}
{\sf E}(Y_{n+1} \vert V_1, V_2, {\cal{D}}) = {\sf E}(Y_{n+1} \vert V_1, {\cal{D}})
\end{equation*}
and
\begin{equation*}
{\sf Var}(Y_{n+1}  \vert V_1, V_2, {\cal{D}})
=
{\sf Var}(Y_{n+1} \vert V_1, {\cal{D}}) .
\nonumber
\end{equation*}
Since ${\sf E}(Y_{n+1} \vert V_1, V_2, {\cal{D}})$ is independent of $V_2$,
\begin{align}
  ~&{\sf E}_{V_1}{\sf Var}_{V_2}E(Y_{n+1} \vert  {\cal{D}}_{n}, V_1, V_2)\nonumber\\
  =&{\sf E}_{V_1 \vert {\cal{D}}} {\sf Var}_{V_2 \vert V_1, {\cal{D}} } E(Y_{n+1} \vert V_1, {\cal{D}} )=0.\nonumber
\end{align}
Hence,
\begin{align}
~&{\sf Var}(Y_{n+1}\vert  {\cal{D}}_{n})\nonumber\\
  =&{\sf E}_{V_1 \vert {\cal{D}} } {\sf E}_{V_2 \vert V_1, {\cal{D}} }{\sf Var} (Y_{n+1} \vert V_1, {\cal{D}})\nonumber\\
 &\quad +{\sf Var}_{V_1 \vert {\cal{D}}} {\sf E}_{V_2 \vert V_1, {\cal{D}} } {\sf E}(Y_{n+1} \vert V_1, {\cal{D}} )\nonumber\\
=&{\sf E}_{V_1 \vert {\cal{D}} } {\sf Var} (Y_{n+1} \vert V_1, {\cal{D}})+Var_{V_1 \vert {\cal{D}}} {\sf E}(Y_{n+1} \vert V_1, {\cal{D}} ). \nonumber
\end{align}
\end{proof}
\end{thm}
In the context of Fig. \ref{fig:horizontal}, 
Theorem \ref{thm:struct} shows that 
conditioning only on the parameter nearest the data is enough.  That is,
given the parameters, the hyperparameters don't matter.  In fact,
the proof of Theorem \ref{thm:struct} shows that it is enough 
for $\cind{Y_{n+1}}{V_2}{\left(V_1,\cal{D}\right)}$ to hold so
Assumption \ref{asmp:condInd} is sufficient, but not necessary. 
On the other hand, Assumption \ref{asmp:condInd} holds in HM,
but not, for instance, in a hidden Markov model.

Now, in the special case of normality, we can show explicitly 
that as you go up a hierarchy satisfying Assumption \ref{asmp:condInd}
the conditional variances increase, at least for the normal.
This means that as the unknown quantity gets further and further from
the data, the data say less and less about it.  We have the following.

\begin{thm}
\label{thm:var} 
Suppose the conditional independence relation in \eqref{condInd} holds. 
Assume that the variables are jointly multivariate Gaussian. Then for 
any choice of the parameters:
  \begin{enumerate}
  \item ${\sf Var}(Y_{n+1}\mid V_2,\mathcal{D}) \ge {\sf Var}(Y_{n+1}\mid V_1,\cal{D})$.
    \item ${\sf Var}_{V_1 \vert {\cal{D}}} [{\sf E} (Y_{n+1} \vert V_1, {\cal{}})]\ge{\sf Var}_{V_2 \vert {\cal{D}}} [{\sf E} (Y_{n+1} \vert V_2, {\cal{D}})]$
    \end{enumerate}
  \begin{proof} 
  Suppose $\Sigma$ is the variance-covariance matrix of all variables and we treat $\cal{D}$ as a multivariate component.  We start with the first clause.
  
\noindent
    Using the formula of conditional covariance for multivariate normals, we get:
    \[
    \sigma_{Y_{n+1}V_2\mid V_1}=\sigma_{Y_{n+1}V_2}-\frac{\sigma_{Y_{n+1}V_1}\sigma_{V_1V_2}}{\sigma_{V_2V_2}}.
        \]
        Now if $\cind{Y_{n+1}}{V_2}{V_1}$, $\sigma_{Y_{n+1}V_2\mid V_1}=0$, so we get
        \begin{equation}\label{eq:cond}
          \sigma_{Y_{n+1}V_2}=\frac{\sigma_{Y_{n+1}V_1}\sigma_{V_1V_2}}{\sigma_{V_2V_2}}.
        \end{equation}
Next, note two identities:
        \begin{align}
          \sigma_{Y_{n+1}Y_{n+1}\mid V_2\mathcal{D}}&=\sigma_{Y_{n+1}Y_{n+1}\mid \mathcal{D}}-\frac{\sigma^2_{Y_{n+1}V_2\mid \mathcal{D}}}{\sigma_{V_2V_2\mid \mathcal{D}}}\nonumber\\
          \sigma_{Y_{n+1}Y_{n+1}\mid V_1\mathcal{D}}&=\sigma_{Y_{n+1}Y_{n+1}\mid \mathcal{D}}-\frac{\sigma^2_{Y_{n+1}V_1\mid \mathcal{D}}}{\sigma_{V_1V_1\mid \mathcal{D}}} .\nonumber
        \end{align}

Now, it is enough to show that
        \begin{align}
        ~&\sigma_{Y_{n+1}Y_{n+1}\mid V_2\mathcal{D}}-\sigma_{Y_{n+1}Y_{n+1}\mid V_1\mathcal{D}}\nonumber\\
        =&\frac{\sigma^2_{Y_{n+1}V_1\mid \mathcal{D}}}{\sigma_{V_1V_1\mid \mathcal{D}}}-\frac{\sigma^2_{Y_{n+1}V_2\mid \mathcal{D}}}{\sigma_{V_2V_2\mid \mathcal{D}}}\ge 0.\nonumber
        \end{align}
We have that
        \[
          \sigma_{Y_{n+1}V_2\mid \mathcal{D}}=\sigma_{Y_{n+1}V_2}-\Sigma_{Y_{n+1}\mathcal{D}}\Sigma^{-1}_{\mathcal{D}\mathcal{D}}\Sigma_{\mathcal{D}V_2}.
          \]
So, from \eqref{condInd} we get $\cind{Y_{n+1}}{V_2}{V_1}$ and 
$\cind{\mathcal{D}}{V_2}{V_1}$.  Similar to the argument above in \eqref{eq:cond} we get $\Sigma_{\mathcal{D}V_2}=\Sigma_{\mathcal{D}V_1}\sigma_{V_1V_2}/\sigma_{V_1V_1}$.  Now, by substitution we get:
          \begin{align}
            \sigma_{Y_{n+1}V_2\mid \mathcal{D}}&=\frac{\sigma_{Y_{n+1}V_1}\sigma_{V_1V_2}}{\sigma_{V_1V_1}}-\Sigma_{Y_{n+1}\mathcal{D}}\Sigma^{-1}_{\mathcal{D}\mathcal{D}}\Sigma_{\mathcal{D}V_1}\frac{\sigma_{V_1V_2}}{\sigma_{V_1V_1}}\nonumber\\
            &=\frac{\sigma_{V_1V_2}}{\sigma_{V_1V_1}}\left\{\sigma_{Y_{n+1}V_1}-\Sigma_{Y_{n+1}\mathcal{D}}\Sigma^{-1}_{\mathcal{D}\mathcal{D}}\Sigma_{\mathcal{D}V_1}\right\}\nonumber\\
            &=\frac{\sigma_{V_1V_2}}{\sigma_{V_1V_1}}\sigma_{Y_{n+1}V_1\mid \mathcal{D}}.\nonumber
          \end{align}
          Furthermore,
          \begin{align}
            \sigma_{V_2V_2\mid \mathcal{D}}&=\sigma_{V_2V_2}-\Sigma_{V_2\mathcal{D}}\Sigma^{-1}_{\mathcal{D}\mathcal{D}}\Sigma_{\mathcal{D}V_2}\nonumber\\
            &=\sigma_{V_2V_2}-\frac{\sigma^2_{V_1V_2}}{\sigma^2_{V_1V_1}}\Sigma_{V_1\mathcal{D}}\Sigma^{-1}_{\mathcal{D}\mathcal{D}}\Sigma_{\mathcal{D}V_1}\nonumber\\
            &=\sigma_{V_2V_2}-\frac{\sigma^2_{V_1V_2}}{\sigma^2_{V_1V_1}}\left(\sigma_{V_1V_1}-\sigma_{V_1V_1\mid \mathcal{D}}\right)\nonumber\\
            &=\sigma_{V_2V_2}-\frac{\sigma^2_{V_1V_2}}{\sigma_{V_1V_1}}+\frac{\sigma^2_{V_1V_2}}{\sigma^2_{V_1V_1}}\sigma_{V_1V_1\mid \mathcal{D}}.\nonumber            
          \end{align}
Re-arranging, we see that
          \begin{equation}\label{eq:end}
            \sigma_{V_2V_2\mid \mathcal{D}}-\frac{\sigma^2_{V_1V_2}}{\sigma^2_{V_1V_1}}\sigma_{V_1V_1\mid \mathcal{D}}=\sigma_{V_2V_2}-\frac{\sigma^2_{V_1V_2}}{\sigma_{V_1V_1}}=\sigma_{V_2V_2\mid V_1}\ge 0.\nonumber
          \end{equation}
So, we get the first clause:
          \begin{align}
            ~&\frac{\sigma^2_{Y_{n+1}V_1\mid \mathcal{D}}}{\sigma_{V_1V_1\mid \mathcal{D}}}-\frac{\sigma^2_{Y_{n+1}V_2\mid \mathcal{D}}}{\sigma_{V_2V_2\mid \mathcal{D}}}\nonumber\\
            =&\frac{\sigma^2_{Y_{n+1}V_1\mid\mathcal{D}}}{\sigma_{V_1V_1\mid\mathcal{D}}}-\frac{\sigma^2_{Y_{n+1}V_1\mid\mathcal{D}}\sigma^2_{V_1V_2}}{\sigma^2_{V_1V_1}\sigma_{V_2V_2\mid\mathcal{D}}}\nonumber\\
              &=\frac{\sigma^2_{Y_{n+1}V_1\mid\mathcal{D}}}{\sigma_{V_1V_1\mid\mathcal{D}}\sigma_{V_2V_2\mid\mathcal{D}}}\left(\sigma_{V_2V_2\mid\mathcal{D}}-\frac{\sigma^2_{V_1V_2}}{\sigma^2_{V_1V_1}}\sigma_{V_1V_1\mid\mathcal{D}}\right)\ge 0.\nonumber
          \end{align}

Clause 2 follows by combining Clause 1 with the observation that
the totals of the terms in the two two-term expansions are equal.
\end{proof}
\end{thm}
We illustrate Theorem \ref{thm:var} with a slightly trivial 
Bayesian hierarchical 
model. 
\begin{exm}
\label{exm:nnn}
Suppose $Y_i$ 
for $i=1, ..., n$
are IID ${\sf N}(\mu, \sigma_0^2)$ where 
$\sigma_0$ is known. Assume $\mu$ is distributed 
as ${\sf N}(\nu, \tau^2_0)$ and $\nu$ is 
distributed as ${\sf N}(a, b^2)$ where both 
$a$ and $b$ are known.
Let
\[
\eta_n = \frac{1}{n/\sigma_0^2 + 1/(\tau_0^2 + b^2)}.
\]
Then, the posterior predictive distribution 
$(Y_{n+1} \vert y^n)$ is
\[
{\sf N}\left(\eta_n \left(\frac{n\bar{y}}{\sigma_0^2} + \frac{a}{\tau_0^2+b^2}\right), \sigma_0^2 + \eta_n\right)
\]
with 
\begin{equation}
\label{3N}
{\sf Var}(Y_{n+1} \vert y^n) = \sigma^2_0 + (n/\sigma_0^2 + 1/(\tau^2_0 + b^2))^{-1}.
\end{equation}
The conditional distribution
$(\mu \vert y^n, \nu)$
is 
$$
{\sf N}\left(\frac{n\bar{Y}/\sigma^2_0 + \nu/\tau_0^2}{n/\sigma_0^2 + 1/\tau_0^2}, \frac{1}{n/\sigma_0^2 + 1/\tau_0^2} \right).
$$
The conditional distribution $(\nu \vert y^n)$ is
$$
{\sf N} \left(\frac{a/b^2 + \bar{y}/(\tau_0^2 + \sigma_0^2/n)}{1/b^2 + 1/(\tau_0^2 + \sigma_0^2/n)}, \frac{1}{1/b^2 + 1/(\tau_0^2 + \sigma_0^2/n)}\right).
$$
Now, it is easy to see that 
\begin{align}
~&{\sf E}_{\nu \vert y^n} 
{\sf E}_{\mu \vert \nu, y^n} 
{\sf Var}(Y_{n+1} \vert y^n, \mu, \nu)
=\sigma_0^2,\nonumber\\
~&{\sf E}_{\nu \vert y^n} 
{\sf Var}_{\mu \vert \nu, y^n} 
{\sf E}(Y_{n+1} \vert y^n, \mu, \nu)\nonumber\\
&=(n/\sigma_0^2 + 1/\tau_0^2)^{-1},\nonumber\\
~&{\sf Var}_{\nu \vert y^n} 
{\sf E}_{\mu \vert \nu, y^n} 
{\sf E}(Y_{n+1} \vert y^n, \mu, \nu)\nonumber\\
&=\frac{(n/\sigma_0^2 + 1/\tau_0^2)^{-2}}
{\tau_0^4(1/b^2 + 1/(\tau_0^2 + \sigma_0^2/n))},
\end{align}
and that these three terms sum to
\eqref{3N}.
As expected,
\begin{align}
{\sf Var}(Y_{n+1} \vert y^n, \mu) = \sigma_0^2
&< \sigma_0^2 + \left( \frac{n}{\sigma^2_0} +
\frac{1}{\tau_0^2} \right)^{-1}\nonumber\\
&={\sf Var}(Y_{n+1} \vert y^n, \nu).
\end{align}\hfill\qed
\end{exm}


As in Theorem \ref{thm:struct}, Assumption \ref{asmp:condInd} is sufficient 
but not necessary for Theorem \ref{thm:var}.  As is evident from the 
proof, the required conditions are $\cind{Y_{n+1}}{V_2}{V_1}$ 
and $\cind{\mc{D}}{V_2}{V_1}$.  It is well-known that 
\eqref{condInd} implies these two conditions, but the converse
does not hold (see \citep{lauritzenBook}).  

Theorem \ref{thm:var} also confirms the fact that a hyperparameter 
has less information about the data and the predicted value than 
a parameter does.  This is similar to the data processing inequality,
see \cite{Cover:Thomas:2006}.   In addition, if the joint density is Gaussian, 
then $\rho^2_{y_{n+1}V_2\vert \mc{D}}\le\rho^2_{y_{n+1}V_1\vert \mc{D}}$, where
$\rho$ is the partial correlation of its subscripts (see \citep{scphd,chaudhuri_2014}).


\subsection{Uncertainty Quantification under posterior independence}

Empirically, often a term in \eqref{LTV2termgen} 
being zero coincides with a term in \eqref{LTV2termgen2} being zero as well.  
We give results showing when this happens.    Consider the following.
\begin{asmp}
\label{asmp:postcond}
  $V_1$ is conditionally independent of $V_2$ given the data $\mc{D}$, i.e., 
  \begin{equation}
  \label{eq:postcond}
    \cind{V_1}{V_2}{\mc{D}}.
    \end{equation}  
\end{asmp}

While constricting, Assumption \ref{asmp:postcond} is satisfied by a large number
of parametric families.   An obvious example is when two experiments are
combined, e.g.,  
$Y_i \sim {\sf Poisson}(\lambda_i)$ and $\lambda_i \sim {\sf Gamma}(a_i, b_i)$
for $i=1, 2$.  This easily extends to many outcomes of $Y$. Another class
of examples is
exponential families whose sufficient statistics split additively 
across parameters equipped with conjugate priors.
A more interesting example is the following.

\begin{exm}
\label{exm:poissontime}
Consider $Y_1 \sim {\sf Poisson}(t \lambda_1)$ and $Y_2 \sim {\sf Poisson}(t \lambda_2 )$,
where $\lambda_i \sim {\sf Gamma} (a_i, b)$, with the $Y_i$'s and the 
$\lambda_i$'s independent. The likelihood factors into separate parts
for $\lambda_1$ and $\lambda_2$, though both factors have $t$.  The
priors are independent, so Assumption \ref{asmp:postcond} holds.
\end{exm}

We have the following implications when a term in a three-term expansion is zero.

\begin{thm}
\label{thm:crsExp}
Suppose Assumption \ref{asmp:postcond} holds. Then we have:
  \begin{enumerate}
  \item If term \eqref{eq:1b} in \eqref{LTV2termgen} is zero, then term 
  \eqref{eq:2c} in \eqref{LTV2termgen2} is zero. That is,
    \begin{align}
      ~&{\sf E}_{V_1 \vert {\cal{D}} } {\sf Var}_{V_2 \vert V_1, {\cal{D}}} [{\sf E} (Y_{n+1} \vert V_1, V_2, {\cal{D}})]=0\nonumber\\
      &\Longrightarrow~ {\sf Var}_{V_2 \vert {\cal{D}}} [{\sf E} (Y_{n+1} \vert V_2, {\cal{D}})]=0.\label{eq:first}
    \end{align}
  \item If term \eqref{eq:2b} in \eqref{LTV2termgen2} is zero, then term \eqref{eq:1c} in \eqref{LTV2termgen} is zero. That is,
    \begin{align}
      ~&{\sf E}_{V_2 \vert {\cal{D}} } {\sf Var}_{V_1 \vert V_2, {\cal{D}}} [{\sf E} (Y_{n+1} \vert V_2, V_1, {\cal{D}})]=0\nonumber\\
      &\Longrightarrow~{\sf Var}_{V_1 \vert {\cal{D}}} [{\sf E} (Y_{n+1} \vert V_1, {\cal{D}})]=0.\label{eq:second}
      \end{align}
  \end{enumerate}
  \begin{proof}
We will only prove Clause I. The proof of Clause II is similar.

Since $\cind{V_1}{V_2}{\mc{D}}$, we have
    \begin{align}
      ~&{\sf Var}_{V_2 \vert V_1, {\cal{D}}} ({\sf E} (Y_{n+1} \vert V_1, V_2, {\cal{D}}))\nonumber\\
      =&{\sf Var}_{V_2 \vert {\cal{D}}} ({\sf E} (Y_{n+1} \vert V_1, V_2, {\cal{D}}))=0.\nonumber
    \end{align}
The LHS of \eqref{eq:first} is zero, so it follows that
    \begin{align}
      &{\sf Var}_{V_2 \vert V_1, {\cal{D}}} [{\sf E} (Y_{n+1} \vert V_1, V_2, {\cal{D}})]=0  \nonumber\\
      \Longrightarrow&{\sf E} (Y_{n+1} \vert V_1, V_2, {\cal{D}})\mbox{~is a constant in terms of $V_2$,}\nonumber\\
      &\hfill\mbox{for $V_1$ and $\mc{D}$.}\nonumber\\
      \Longrightarrow&\frac{\partial}{\partial V_2}{\sf E} (Y_{n+1} \vert V_1, V_2, {\cal{D}})=0 \mbox{  $\forall$ $V_1$, $V_2$ and $\mc{D}$.}\nonumber
      \end{align}
    Now, using the Leibnitz rule for all $V_1$, $V_2$, and $\mc{D}$:
    \begin{align}
   & \frac{\partial}{\partial V_2} E(Y_{n+1} \vert V_2, \mc{D}) =
      ~\frac{\partial}{\partial V_2}{\sf E}_{V_1 \vert {\cal{D}} }{\sf E} (Y_{n+1} \vert V_1, V_2, {\cal{D}})\nonumber\\
      &=\frac{\partial}{\partial V_2}\int {\sf E} (Y_{n+1} \vert V_1, V_2, {\cal{D}})f_{V_1 \vert {\cal{D}} } dV_1\nonumber\\
      &=\int\frac{\partial}{\partial V_2} {\sf E} (Y_{n+1} \vert V_1, V_2, {\cal{D}})f_{V_1 \vert {\cal{D}} } dV_1=0.\nonumber
    \end{align}
    That is, ${\sf E}_{V_1 \vert {\cal{D}} }{\sf E} (Y_{n+1} \vert V_1, V_2, {\cal{D}})={\sf E} [Y_{n+1} \vert V_2, {\cal{D}}]$ is a constant in terms of $V_2$, for all $V_1$ and $\mc{D}$.  This implies that ${\sf Var}_{V_2 \vert {\cal{D}}} [{\sf E} (Y_{n+1} \vert V_2, {\cal{D}})]=0$.
    \end{proof}
\end{thm}
Under Assumption \ref{asmp:postcond}, to satisfy the condition on the left
of \eqref{eq:first}, it is enough for ${\sf E} (Y_{n+1} \vert V_1, V_2, {\cal{D}})$ 
to be independent of $V_2$ for all $V_1$ and $\mc{D}$.  The conditional independence $\cind{Y_{n+1}}{V_2}{(V_1,\mc{D})}$ is not required.  
In fact, ${\sf Var}(Y_{n+1} \vert V_1, V_2, {\cal{D}})$ may still depend on $V_2$.  
The analogous statements hold for the condition on the left of \eqref{eq:second}.

Assumption \ref{asmp:postcond} is sufficient but not necessary for 
Theorem \ref{thm:crsExp} to hold.  The next two examples show that
i) without Assumption \ref{asmp:postcond} Theorem \ref{thm:crsExp}
need not hold, and
ii) the conclusions of Theorem \ref{thm:crsExp} can hold even
when Assumption \ref{asmp:postcond} does not.

\begin{exm}
\label{exm:norm1} 
Let $Y_1$, $Y_2$, $\ldots$, $Y_n$ be IID ${\sf N}\left(\mu,1/\lambda^2\right)$ with 
$\mu\sim {\sf N} \left(\mu_0,1/\lambda^2_0\right)$ with known $\mu_0$ and 
$\lambda_0$ and $\lambda^2\sim {\sf Gamma}\left(\alpha_0,\beta_0\right)$ 
with known $\alpha_0$ and $\beta_0$.
Set $V_1=\mu$, $V_2=\lambda^2$, and $\mc{D}=\left\{Y_1,Y_2,\ldots,Y_n\right\}$.  
Then we can show that 
$$
\mu\vert \lambda^2,\mc{D} \sim {\sf N} \left(\frac{\lambda^2\sum^n_{i=1}Y_i+\lambda^2_0\mu_0}{n\lambda^2+\lambda^2_0},\frac{1}{n\lambda^2+\lambda^2_0}\right).
$$  
Since this distribution depends on $\lambda^2$, $\mu\not\upmodels\lambda^2\vert\mc{D}$ and Assumption \ref{asmp:postcond} does not hold.

We see that Theorem \ref{thm:crsExp} does not hold either.
Even though ${\sf E}[Y_{n+1}\vert\mu,\lambda^2,\mc{D}]=\mu$, is free of $\lambda^2$,
i.e.,
  \[
   {\sf E}_{\mu\vert {\cal{D}} } {\sf Var}_{\lambda^2 \vert \mu, {\cal{D}}} [{\sf E} (Y_{n+1} \vert \mu, \lambda^2, {\cal{D}})]={\sf E}_{\mu \vert {\cal{D}} } {\sf Var}_{\lambda^2 \vert \mu, {\cal{D}}} [\mu]=0
   \]
   we also have that
   \begin{align}
     ~&{\sf Var}_{\lambda^2 \vert {\cal{D}}} ({\sf E} (Y_{n+1} \vert \lambda^2, {\cal{D}}))\nonumber\\
     =&{\sf Var}_{\lambda^2 \vert {\cal{D}}} {\sf E}_{\mu\vert \cal{D},\lambda^2 }({\sf E} (Y_{n+1} \vert \mu,\lambda^2, {\cal{D}}))\\
     =&{\sf Var}_{\lambda^2 \vert {\cal{D}}}{\sf E}_{\mu\vert \cal{D},\lambda^2 }(\mu)\nonumber\\
     =&{\sf Var}_{\lambda^2 \vert {\cal{D}}}\left(\frac{\lambda^2\sum^n_{i=1}Y_i+\lambda^2_0\mu_0}{n\lambda^2+\lambda^2_0}\right)>0.
   \end{align}
   The last inequality holds because ${\sf E}(\mu\vert \mc{D},\lambda^2)$ 
   is not free of $\lambda^2$.\hfill\qed
  \end{exm}

\begin{exm}
\label{exm:norm2} 
Suppose $\mu\sim {\sf N}\left(\mu_0,1/(\kappa_0\lambda^2)\right)$ is used in
Example \ref{exm:norm1} so that Assumption \ref{asmp:postcond} still does not hold.
Now,
\[
\mu\vert\lambda^2,\mc{D} \sim {\sf N}\left(\frac{(\sum^n_{i=1}Y_i+\kappa_0\mu_0)}{(n+\kappa_0)},\frac{1}{(n+\kappa_0)\lambda^{2}}\right)
\]
and $\mu$ remains
dependent on $\lambda^2$ (given $\mc{D}$). 
We also see that $E[Y_{n+1}\vert\mu,\lambda^2,\mc{D}]=\mu$ remains free of 
$\lambda^2$
and this gives
  \begin{align}
    {\sf E}_{\mu\vert {\cal{D}} } {\sf Var}_{\lambda^2 \vert \mu, {\cal{D}}} [{\sf E} (Y_{n+1} \vert \mu, \lambda^2, {\cal{D}})]&={\sf E}_{\mu \vert {\cal{D}} }{\sf Var}_{\lambda^2 \vert \mu, {\cal{D}}} (\mu)\nonumber\\
    &=0.\nonumber
   \end{align}
In addition, unlike Example \ref{exm:norm1}, we find
   \begin{align}
     ~&{\sf Var}_{\lambda^2 \vert {\cal{D}}}[{\sf E} (Y_{n+1} \vert \lambda^2, {\cal{D}})]={\sf Var}_{\lambda^2 \vert {\cal{D}}}{\sf E}_{\mu\vert \cal{D},\lambda^2 }(\mu)\nonumber\\
     =&{\sf Var}_{\lambda^2 \vert {\cal{D}}}\left(\frac{\sum^n_{i=1}Y_i+\kappa_0\mu_0}{n+\kappa_0}\right)=0.\nonumber
   \end{align}
That is, the conclusions of Theorem \ref{thm:crsExp} are satisfied.\hfill\qed
  \end{exm}

Even though Assumption \ref{asmp:postcond} is not necessary, from Examples \ref{exm:norm1} and \ref{exm:norm2} it is evident that the necessary condition is difficult to specify, might depend on the parametrisation, and be
quite hard to verify in practice.   Indeed, to make term \eqref{eq:1c}
zero, without $\cind{V_1}{V_2}{\mc{D}}$,
we effectively need ${\sf E}\left[Y_{n+1}\vert V_1,V_2,\mc{D}\right]=\mu(V_1,\mc{D})$
and ${\sf E}_{V_1\vert V_2,\mc{D}}[\mu(V_1,\mc{D})]$ to be free of $V_2$.
Example \ref{exm:norm1} clearly shows that, depending on the parametrisation, such 
conditions may or may not hold. Moreover, they may be hard to verify and
interpret, whereas the sufficient Assumption \ref{asmp:postcond} can be 
relatively easily verified and interpreted, e.g., in Example \ref{exm:hmm}.

It is easy to see that under Assumption \eqref{asmp:postcond}, neither \eqref{eq:first} nor \eqref{eq:second} leads to an interpretable two-term 
expansion of the predictive variance.  For instance, if we have
${\sf E}_{V_1 \vert {\cal{D}} } {\sf Var}_{V_2 \vert V_1, {\cal{D}}} [{\sf E} (Y_{n+1} \vert V_1, V_2, {\cal{D}})]=0$
then we get
\begin{align}
  {\sf Var}( Y_{n+1} \vert {\cal{D}})=&{\sf E}_{V_1,V_2 \vert {\cal{D}}} ( {\sf Var} ( Y_{n+1}\vert V_2, V_1, {\cal{D}}))\nonumber\\
  +&{\sf Var}_{V_1 \vert {\cal{D}}} ({\sf E} (Y_{n+1} \vert V_1, {\cal{D}})),
\end{align}
in which ${\sf Var} ( Y_{n+1}\vert V_2, V_1, {\cal{D}})$ is not in
general independent of $V_2$, similarly if $V_1$ and $V_2$ are interchanged.

The converses in Theorem \ref{thm:crsExp} do not hold either without 
extra assumptions.  We have the
following result, which is symmetric in $V_1$ and $V_2$.

\begin{thm}
\label{thm:rev} 
Suppose Assumption \ref{asmp:postcond} holds.  Then:
  \begin{enumerate}
  \item If, additionally, for all $V_1$, $V_2$, $\mc{D}$, 
  we have that $\cind{Y_{n+1}}{V_2}{(V_1,\mc{D})}$, then
    \[
      {\sf E}_{V_1 \vert {\cal{D}} } {\sf Var}_{V_2 \vert V_1, {\cal{D}}} [{\sf E} (Y_{n+1} \vert V_1, V_2, {\cal{D}})]=0
      \]
      and
      \[
      {\sf Var}_{V_2 \vert {\cal{D}}} [{\sf E} (Y_{n+1} \vert V_2, {\cal{D}})]=0.
    \]
  
      \item In this case, the three-term expansion using the LTV on $V_1$ 
first and $V_2$ second reduces to a two-term expansion.
  \end{enumerate}
  
  \begin{proof} 
  For Clause 1, if $\cind{Y_{n+1}}{V_2}{(V_1,\mc{D})}$ for all $V_1$, $V_2$, $\mc{D}$, then it trivially follows that ${\sf E} [Y_{n+1} \vert V_1, V_2, {\cal{D}}]={\sf E} [Y_{n+1} \vert V_1, {\cal{D}}]$, which is independent of $V_2$.  That is, we have both
\begin{align}
  ~&{\sf Var}_{V_2 \vert V_1, {\cal{D}}} [{\sf E} (Y_{n+1} \vert V_1, V_2, {\cal{D}}])\nonumber\\
=&{\sf Var}_{V_2 \vert {\cal{D}}} ({\sf E} (Y_{n+1} \vert V_1, {\cal{D}}))=0
\nonumber
\end{align}
and
\begin{align}
~&{\sf Var}_{V_2 \vert {\cal{D}}} [{\sf E} (Y_{n+1} \vert V_2, {\cal{D}})]\nonumber\\
=&{\sf Var}_{V_2 \vert {\cal{D}}} [{\sf E}_{V_1\vert V_2,\mc{D}}{\sf E} (Y_{n+1} \vert V_1, V_2, {\cal{D}})]\nonumber\\
=&{\sf Var}_{V_2 \vert {\cal{D}}} [{\sf E}_{V_1\vert \mc{D}}{\sf E} (Y_{n+1} \vert V_1, {\cal{D}})]=0.\nonumber
\end{align}
Clause 2 follows from Theorem \ref{thm:struct}.
\end{proof}
\end{thm}

From the proof of Theorems \ref{thm:crsExp} and \ref{thm:rev}, we 
see that we only require 
${\sf E} [Y_{n+1} \vert V_1, V_2, {\cal{D}}]$ to be independent of $V_2$ 
for the equivalence of the three and two-term expansions.  This is the
familiar condition first-order ancillarity, \cite{Lehmann:1998}, p. 41.
For Clause 2, we also require that 
${\sf Var}[Y_{n+1} \vert V_1, V_2, {\cal{D}}]$ be independent 
of $V_2$.  We assume this stronger structural condition 
$\cind{Y_{n+1}}{V_2}{(V_1,\mc{D})}$ because it is easier to
verify.

Note that the assumption that $\cind{V_1}{V_2}{\mc{D}}$ (Assumption
\ref{asmp:postcond}) and $\cind{Y_{n+1}}{V_2}{(V_1,\mc{D})}$ imply that the 
condition $\cind{(Y_{n+1},V_1)}{V_2}{\mc{D}}$ holds for all $V_1$, $V_2$, 
and $\mc{D}$ (see \citep{lauritzenBook}).  For many models, such relationships can be easily 
determined from their description. 
Theorems \ref{thm:crsExp} and \ref{thm:rev} would rarely apply to a 
hierarchical Bayes model. However, the hidden Markov Model in 
Example \ref{exm:hmm} would satisfy all conditions of both theorems. 

\begin{exm} These conditional independence relations are satisfied by the hidden Markov model described in Example \ref{exm:hmm}. By construction, we have $(V_1 \upmodels V_2 \vert {\cal{D}})$ and
$(Y_{n+1} \upmodels V_2 \vert {\cal{D}}, V_1)$.  Now using Theorem \ref{thm:rev} it immediately follows that term \eqref{eq:1b} in \eqref{LTV2termgen} and
term \eqref{eq:2b} in \eqref{LTV2termgen} are zero, and the three-term 
expansions reduce to a two-term expansion involving $V_2$.  

Furthermore, since the predicted values $\hat{Y}_i$ of $Y_i$
that are functions of the $X_1$, $X_2$, $\ldots$, $X_i$, by letting $\widehat{\mc{D}}=\left\{\hat{Y}_1,\hat{Y}_2,\ldots,\hat{Y}_n\right\}$, the conditions $\cind{V_2}{V_1}{\widehat{\mc{D}}}$ and $\cind{Y_{n+1}}{V_1}{(V_2,\widehat{\mc{D})}}$ still hold. That is, Theorem \ref{thm:rev} applies.\qed
\end{exm}

\section{General Expansions of the PPV}
\label{decomposition}

Here, we extend the results of the previous sections to general
multi-term expansions of the posterior predictive variance.  
Without loss of generality, we fix $\mc{V}=\{V_1,V_2,\ldots, V_K\}$ to 
be a specific ordering of the
entries in $\mc{Z}$ and assume that each element in $\mc{Z}$ 
appears in our expansion.  That is, each $Z_i$ is \emph{manifest}, no 
element is \emph{latent}.  

A trivial but condensed two-term expansion of the PPV given $\mc{V}$ and $\mc{D}$ is
\begin{align}
\label{condensed_var}
~&\nonumber {\sf Var}(Y_{n+1} \vert  {\cal{D}}_n)(\mathcal{V})\nonumber\\ 
 =& {\sf E}_{(V_1,\ldots, V_K)} {\sf Var}(Y_{n+1} \vert {\cal{D}}_n, V_1, \ldots, V_K) \\ 
& + {\sf Var}_{(V_1,\ldots, V_K)}{\sf E}(Y_{n+1} \vert  {\cal{D}}_n, V_1, \ldots, V_K).
\end{align}
More interesting is the $K+1$-term expansion of posterior predictive variance ${\sf Var}[Y_{n+1}\vert \mc{D}]$ w.r.t. $\mc{V}$.  This is given by:
\begin{subequations}\label{Conditional_Var_sum}
\begin{align}
~&{\sf Var}(Y_{n+1} \vert  {\cal{D}}_n)(\mc{V})\nonumber\\ 
=&{\sf E}_{(V_1,\ldots, V_K)\vert \mc{D}} {\sf Var}(Y_{n+1} \vert  {\cal{D}}_n, V_1, \ldots, V_K)
\label{eq:g1}\\  
+&\sum_{k=K}^{2} {\sf E}_{(V_{1}, \ldots, V_{k-1} )\vert\mc{D}} {\sf Var}_{V_k\vert\mc{D},V_{1},\ldots,V_{k-1}}\left[\right.\nonumber\\
&\left.\quad\quad\quad{\sf E}(Y_{n+1} \vert  {\cal{D}}_n, V_1, \ldots, V_{k})\right]
\label{eq:g2}\\
+&{\sf Var}_{V_1\vert\mc{D}}{\sf E}(Y_{n+1} \vert  {\cal{D}}_n, V_1).
\label{eq:g3}
\end{align}
\end{subequations} 
The terms in \eqref{eq:g2} are added in decreasing order of $k$ for 
notational convenience; the motivation for this will become clear in the sequel.

The inner expectation in the $k$-th in summands in \eqref{eq:g2} and 
\eqref{eq:g3} can be iterated as 
a sequence of conditional expectations:
\begin{align}
 ~&{\sf E}(Y_{n+1} \vert  {\cal{D}}_n, V_1, \ldots, V_{k})\nonumber\\
  =&{\sf E}_{V_{k+1}\vert \mc{D}, V_{1},\ldots,V_{k}}\cdots{\sf E}_{V_{K}\vert \mc{D}, V_{1},\ldots,V_{K-1}}\left[\right.\nonumber\\
  &\left.\quad{\sf E}(Y_{n+1} \vert  {\cal{D}}_n, V_1, \ldots, V_{K})\right]\nonumber
\end{align}
and
\begin{align}
  ~&{\sf E}(Y_{n+1} \vert  {\cal{D}}_n, V_1)\nonumber\\
  =&{\sf E}_{V_{2}\vert \mc{D}, V_{1}}\cdots E_{V_{K}\vert \mc{D}, V_{1},\ldots,V_{K-1}} {\sf E}(Y_{n+1} \vert  {\cal{D}}_n, V_1, \ldots, V_{K}),\nonumber
\end{align}
respectively.
That is, even though all $V_k$'s appear in the PPV, only $V_1$, $V_2$, $\ldots$, $V_{k}$ appear in the $k$-th term of the sum in \eqref{eq:g2}; the rest of the
elements are latent. 

The expansion in \eqref{Conditional_Var_sum} depends on the choice of $\mc{V}$, i.e., the specific permutation of the elements of $\mc{Z}$.  The results of 
Section \ref{sec:quant} extend to this general expansion whose terms quatify the uncertainty associated with the PPV.  

\subsection{Uncertainty Quantification under Structural Conditions}
\label{UQstructural}

Suppose that the set of variables $\mc{V}$ can be split into 
$\mc{V}_1=\{V_1,V_2,\ldots,V_m\}$ and
$\mc{V}_2=\{V_{m+1},V_{m+2},\ldots,V_K\}$.  We
extend Assumption \ref{asmp:condInd} to the general $K$ 
and $m$ setting.
\begin{asmp}
extends\label{asmp:condindext}
$Y_{n+1}$ and $\cal{D}$ are conditionally 
independent of ${\mc{V}_2}$ given ${\mc{V}_1}$, i.e.,
\begin{equation}
\label{condIndmath}
\cind{\left(Y_{n+1},\mc{D}\right)}{\mc{V}_2}{\mc{V}_1}
\end{equation}
\end{asmp}
Under this assumption, Theorem \ref{thm:struct} admits 
a straightforward extension to the reduction of 
a $(K+1)$-term expansion to an $(M+1)$-term expansion 
given $\mc{V}_1$ and $\mc{D}$.

\begin{thm}
\label{thm:gStCond}
Under Assumption \ref{asmp:condindext} we have 
\[
{\sf Var}[Y_{n+1} \vert {\cal{D}}_n](\mc{V})={\sf Var}[Y_{n+1} \vert {\cal{D}}_n](\mc{V}_1).
\]
\begin{proof}  
Since the conditional independence implies 
$\cind{Y_{n+1}}{\mc{V}_2}{(\mc{D},\mc{V}_1)}$, 
we have 
\begin{equation*}
  {\sf E}\left(Y_{n+1}\vert\mc{D},\mc{V}\right)={\sf E}\left(Y_{n+1}\vert\mc{D},\mc{V}_1\right)
  \end{equation*}
and
\begin{equation*}
{\sf Var}\left(Y_{n+1}\vert\mc{D},\mc{V}\right)={\sf Var}\left(Y_{n+1}\vert\mc{D},\mc{V}_1\right).
\nonumber
\end{equation*}
Hence, in \eqref{eq:g1}, 
\begin{align}
  ~&{\sf E}_{\mc{V}\vert \mc{D}} {\sf Var}(Y_{n+1} \vert  {\cal{D}}_n, \mc{V})\nonumber\\
  =&{\sf E}_{\mc{V}\vert \mc{D}} {\sf Var}(Y_{n+1} \vert  {\cal{D}}_n, \mc{V}_1)={\sf E}_{\mc{V}_1\vert \mc{D}} {\sf Var}(Y_{n+1} \vert  {\cal{D}}_n,\mc{V}_1).\nonumber
\end{align}
Similarly, it follows that:
\begin{align}
  ~&{\sf E}[Y_{n+1} \vert  {\cal{D}}_n, V_1]\nonumber\\
  =&{\sf E}_{V_{2}\vert \mc{D}, V_{1}}\cdots {\sf E}_{V_{K}\vert \mc{D}, V_{1},\ldots,V_{K-1}} {\sf E}[Y_{n+1} \vert  {\cal{D}}_n, \mc{V}]\nonumber\\
  =&{\sf E}_{V_{2}\vert \mc{D}, V_{1}}\cdots{\sf E}_{V_{M}\vert \mc{D}, V_{1},\ldots,V_{M-1}} {\sf E}_{V_{M+1}\vert \mc{D}, \mc{V}_{1}}\cdots\nonumber\\
  &\quad \cdots {\sf E}_{V_{K}\vert \mc{D}, \mc{V}_1,V_{m+1},\cdots,V_{K-1}} {\sf E}(Y_{n+1} \vert  {\cal{D}}_n, \mc{V}_1)\nonumber\\
=&{\sf E}_{V_{2}\vert \mc{D}, V_{1}}\cdots {\sf E}_{V_{M}\vert \mc{D}, V_{1},\ldots,V_{M-1}}{\sf E}(Y_{n+1} \vert  {\cal{D}}_n, \mc{V}_1).
\nonumber
\end{align}
The last equality holds because 
${\sf E}(Y_{n+1} \vert  {\cal{D}}_n, \mc{V}_1)$ is free of elements in $\mc{V}_2$.

Using the above argument we also see that $k>m$ implies the 
summands in \eqref{eq:g2} translate to:
\begin{align}
~&{\sf E}_{(V_{1}, \ldots, V_{k-1} )\vert\mc{D}} {\sf Var}_{V_k\vert\mc{D},V_{1},\ldots,V_{k-1}} {\sf E}(Y_{n+1} \vert  {\cal{D}}_n, V_1, \ldots, V_{k})\nonumber\\
  =&{\sf E}_{(V_{1}, \ldots, V_{k-1} )\vert\mc{D}} {\sf Var}_{V_k\vert\mc{D},\mc{V}_{1},V_{m+1}\ldots,V_{k-1}} {\sf E}(Y_{n+1} \vert  {\cal{D}}_n, \mc{V}_1)\nonumber\\
  =&0,\nonumber
    \end{align}
since ${\sf E}[Y_{n+1} \vert  {\cal{D}}_n, \mc{V}_1]$ is free of $V_k$.

Now, collecting all the terms, we get:
\begin{align}
~&{\sf Var}(Y_{n+1} \vert  {\cal{D}}_n)(\mc{V})\nonumber\\
=&{\sf E}_{\mc{V}_1\vert \mc{D}} {\sf Var}(Y_{n+1} \vert  {\cal{D}}_n, \mc{V}_1)\nonumber\\
+&\sum_{k=M}^{2} E_{(V_{1}, \ldots, V_{k-1})\vert\mc{D}} {\sf Var}_{V_k\vert\mc{D},V_{1},\ldots,V_{k-1}}{\sf E}(Y_{n+1} \vert  {\cal{D}}_n, \mc{V}_1)\nonumber\\
& + {\sf Var}_{V_1\vert\mc{D}} {\sf E}_{V_{2}\vert \mc{D}, V_{1}}\cdots {\sf E}_{V_{M}\vert \mc{D}, V_{1},\ldots,V_{M-1}}{\sf E}(Y_{n+1} \vert  {\cal{D}}_n, \mc{V}_1)\nonumber\\
=&{\sf Var}(Y_{n+1} \vert {\cal{D}}_n)(\mc{V}_1).\nonumber
\end{align}
\end{proof}
\end{thm}

\subsection{Uncertainty Quantification under Posterior Independence}

We now examine how the overall uncertainty as measured by the PPV
spreads over the terms in different permutations of the elements 
of $\mc{Z}$. Let $\pi$ be a permutation of $\{1, \ldots, K\}$ and
$\mc{V}^{\pi}$ denote the corresponding permutation of the 
elements in $\mc{V}$. 
Our goal is to identify conditions under which a given term in one
expansion of $Var[Y_{n+1}\vert\mc{D}](\mc{V})$ is zero implies
that a term in another expansion of $Var[Y_{n+1}\vert\mc{D}](\mc{V}^{\pi})$
is also zero.

Recall $\mc{V}$ is the fixed ordering of the elements in $\mc{Z}$ and
denote $I=\{1,2,\ldots,K\}$.  We write the posterior predictive variance in \eqref{Conditional_Var_sum} as:
\begin{equation}\label{eq:Ts}
    {\sf Var}(Y_{n+1}\vert\mc{D})(\mc{V})=T^I_0+T^I_K+\cdots+T^I_1,
\end{equation}
where
\begin{align}
    T^I_0&={\sf E}_{\mc{V}\vert\mc{D}}{\sf Var}\left(Y_{n+1}\vert\mc{D},\mc{V}\right),\nonumber\\
    T^I_k&={\sf E}_{(V_{1}, \ldots, V_{k-1} )\vert\mc{D}} {\sf Var}_{V_k\vert\mc{D},V_{1},\ldots,V_{k-1}}\left[\right.\nonumber\\
      &\left.\quad{\sf E}(Y_{n+1} \vert  {\cal{D}}_n, V_1, \ldots, V_{k})\right],\nonumber\\
    &\hspace{.5\linewidth}~\mbox{    for $k=2$, $3$, $\ldots$, $K$},\nonumber\\
    T^I_1&={\sf Var}_{V_1\vert\mc{D}}{\sf E}(Y_{n+1} \vert  {\cal{D}}_n, V_1).\nonumber
    \end{align}
Note that $T^I_0$ has the conditional variance of $Y_{n+1}$
while $T^I_k$ has the conditional variance of $V_k$ when for $k\ge1$.   
Write $\mc{V}^{\pi}=\left\{V^{\pi}_1, V^{\pi}_2,\ldots,V^{\pi}_K\right\}$ for 
the permuted elements of $\mc{Z}$, where $V^{\pi}_j=V_i$ if and only if $\pi(i)=j$. 
Similar to \eqref{eq:Ts} we write
\begin{equation}\label{eq:Tps}
{\sf Var}(Y_{n+1}\vert\mc{D})(\mc{V}^{\pi})=T^{\pi}_0+T^{\pi}_K+\cdots+T^{\pi}_1,
\end{equation}
where the term $T^{\pi}_k$ involves computing the variance 
with respect $V^{\pi}_k$, for $k\ge 1$.  

To generalize Theorem \ref{thm:crsExp} to a $(K+1)$-term expansion, we start with 
the following.  Denote by 
\[
\neg i=\left\{K,K-1,\ldots,i+1,i-1,\ldots,1\right\}.  
\]

The following Lemma will be used in Theorem \ref{thm:genExp}.

\begin{Lem}
\label{condind}
Suppose $\cind{V_i}{V_{\neg i}}{\mc{D}}$.  Then, $\forall u=i+1,\ldots K$
\[
\cind{V_i}{V_{u}}{\left(\{\mc{D},V_1,V_2,\ldots,V_{u-1}\}\setminus\{V_i\}\right)}.
\]

\begin{proof}
Consider $u=K$.  The result follows by factoring $V_{\neg i}$ into two factors, one for $K$ and one for $1,\ldots , i-1, i+1, \ldots , K-1$.  Repeating 
this for $K-1$ etc. gives the result.
\end{proof}
\end{Lem}

Now we relate the expansions for ${\sf Var}(Y_{n+1}\vert\mc{D})(\mc{V})$ 
and ${\sf Var}(Y_{n+1}\vert\mc{D})(\mc{V}^{\pi})$.
The first statement says that iterating an expectation that occurs before
a variance does not change a term.  It is understood that the same applies
to the expectations taken after the variance.  The second statement
gives a condition under which if a term is zero, we can,
loosely, say that terms with the variance moved to the left are also zero.

\begin{thm}
\label{thm:genExp}
For the $K+1$ term expansions in
\eqref{eq:Ts} and \eqref{eq:Tps} we have:
\begin{enumerate}

\item For all $i = 1, \ldots, K$, if $\pi(i)=i$, and $\pi\left(\left\{i-1,\ldots,1\right\}\right)=\left\{i-1,\ldots,1\right\}$, then
$T^I_i=0$ $\Longleftrightarrow$ $T^{\pi}_i=0$.

\item Under the conditions of Lemma \ref{condind},
if $\pi(i)=j$, and 
$\left\{j-1,\ldots,1\right\}\subset\pi\left(\left\{i-1,\ldots,1\right\}\right)$, 
we have that
$T^I_i=0$ $\Longrightarrow$ $T^{\pi}_j=0$.

\end{enumerate}
\begin{proof}
For the first statement, note that 
\begin{align}
  ~&T^{I}_i=0\nonumber\\
  &\Longleftrightarrow {\sf Var}_{V_i\vert \mc {D}, V_1, V_2,\ldots, V_{i-1}}{\sf E}[Y_{n+1}\vert\mc{D}, V_1, V_2, \ldots, V_i]=0\nonumber\\
&\Longleftrightarrow {\sf E}[Y_{n+1}\vert\mc{D}, V_1, V_2, \ldots, V_i]\mbox{ is free of $V_i$}\nonumber\\
&\Longleftrightarrow {\sf E}[Y_{n+1}\vert\mc{D}, V^{\pi}_1, V^{\pi}_2, \ldots, V^{\pi}_i]\mbox{ is free of $V^{\pi}_i$}\label{eq:perm}\\
&\Longleftrightarrow {\sf Var}_{V^{\pi}_i\vert \mc {D}, V^{\pi}_1, V^{\pi}_2,\ldots, V^{\pi}_{i-1}}{\sf E}[Y_{n+1}\vert\mc{D}, V^{\pi}_1, V^{\pi}_2, \ldots, V^{\pi}_i]=0\nonumber\\
&\Longleftrightarrow T^{\pi}_i=0.\nonumber
\end{align}

For the second statement, because variance is non-negative we have
\begin{align}
  ~&T^{I}_i=0\nonumber\\
  &\Longleftrightarrow {\sf Var}_{V_i\vert \mc {D}, V_1, V_2,\ldots, V_{i-1}}{\sf E}(Y_{n+1}\vert\mc{D}, V_1, V_2, \ldots, V_i)=0\nonumber\\
&\Longleftrightarrow {\sf E}(Y_{n+1}\vert\mc{D}, V_1, V_2, \ldots, V_i)\mbox{ is free of $V_i$}\nonumber\\
&\Longleftrightarrow \frac{\partial}{\partial V_i}{\sf E}(Y_{n+1}\vert\mc{D}, V_1, V_2, \ldots, V_i)=0 \mbox{ for all $\mc{D}$},
\label{eq:prm}
\end{align}
with mild abuse of notation. 

Next, since $j \leq i$ and $\pi^{-1}(1, \ldots, j)\subseteq \{1, \ldots, i\}$
we fill out $\pi^{-1}(1, \ldots, j)$ to have cardinality $i$ by
writing $\{ 1, \ldots, i\} = \pi^{-1}(\{1, \ldots, j, \ell_1, \ldots \ell_{i-j}\})$.
For ease of notation we express this as $\ell_{i-j} = j^\prime,
\ell_{i-j-1} = j^\prime -1, \ldots , \ell_1 = j+1$, i.e.,
\[
\pi\left(\left\{i-1,\ldots,1\right\}\right)=\left\{j^{\prime},j^{\prime}-1,\ldots,j+1,j-1,\ldots,1\right\}.
\]
From the discussion above:
\begin{align}
&{\sf E}(Y_{n+1}\vert \mc{D}, V^{\pi}_1,\ldots, V^{\pi}_j)
\nonumber \\
&={\sf E}_{V^{\pi}_{j^{\prime}}\vert \mc{D},V^{\pi}_1,\ldots,V^{\pi}_{j^{\prime}-1}}\cdots E_{V^{\pi}_{j+1}\vert \mc{D}, V^{\pi}_1,\ldots,V^{\pi}_{j}}\nonumber\\
&\quad{\sf E}(Y_{n+1}\vert \mc{D}, V^{\pi}_1,\ldots, V^{\pi}_{j-1}, V^{\pi}_{j}, V^{\pi}_{j+1},\ldots,V^{\pi}_{j^{\prime}}).
\label{eq:perm2}
\end{align}


Since $V^{\pi}_j=V_i$, and rest of $\mc{V}^{\pi}$ is a relabelling of $\mc{V}$, the condition in Lemma \ref{condind} translates to $\cind{V^{\pi}_j}{V^{\pi}_{\neg j}}{\mc{D}}$.  
So, by Lemma \ref{condind},
each expectation other than
\[
  {\sf E}(Y_{n+1}\vert \mc{D}, V^{\pi}_1,\ldots, V^{\pi}_{j-1}, V^{\pi}_{j}, V^{\pi}_{j+1},\ldots,V^{\pi}_{j^{\prime}})
  \]
  in \eqref{eq:perm2} is free of $V^{\pi}_j$.  Now, the dominated derivative theorem
gives
\begin{align}
~&\frac{\partial}{\partial V^{\pi}_j}{\sf E}(Y_{n+1}\vert \mc{D}, V^{\pi}_1,\ldots, V^{\pi}_j) 
\nonumber\\
=&
{\sf E}_{V^{\pi}_{j^{\prime}}\vert \mc{D}, V^{\pi}_1,\ldots,V^{\pi}_{j^{\prime}-1}}
\cdots {\sf E}_{V^{\pi}_{j+1}\vert \mc{D}, V^{\pi}_1,\ldots,V^{\pi}_{j-1}}\nonumber\\
&\quad\frac{\partial}{\partial V^{\pi}_{j}}{\sf E}(Y_{n+1}\vert \mc{D}, V^{\pi}_1,\ldots, V^{\pi}_{j-1}, V^{\pi}_{j}, V^{\pi}_{j+1},\ldots,V^{\pi}_{j^{\prime}}). \nonumber
\end{align}
Now, \eqref{eq:prm} implies that:
\begin{align}
~&\frac{\partial}{\partial V^{\pi}_j}{\sf E}(Y_{n+1}\vert \mc{D}, V^{\pi}_1,\ldots, V^{\pi}_{j-1}, V^{\pi}_{j}, V^{\pi}_{j+1},\ldots,V^{\pi}_{j^{\prime}})\nonumber\\
  =&\frac{\partial}{\partial V_i}{\sf E}(Y_{n+1}\vert\mc{D}, V_1, V_2, \ldots, V_i)=0\nonumber\\
  &\hspace{.5\linewidth}\mbox{ for all $\mc{D}$, $V_1$, $V_2$, $\ldots$, $V_i$}.\nonumber
\end{align}
That is:
\begin{align}
~&\frac{\partial}{\partial V^{\pi}_j}{\sf E}(Y_{n+1}\vert \mc{D}, V^{\pi}_1,\ldots, V^{\pi}_j)=0\nonumber\\
\Longleftrightarrow& {\sf E}(Y_{n+1}\vert \mc{D}, V^{\pi}_1,\ldots, V^{\pi}_j)\mbox{ is free of $V^{\pi}_j$}\nonumber\\
\Longleftrightarrow & T^{\pi}_j= {\sf Var}_{V^{\pi}_j\vert\mc{D},V^{\pi}_1,\ldots,V^{\pi}_{j-1}}{\sf E}(Y_{n+1}\vert \mc{D}, V^{\pi}_1,\ldots, V^{\pi}_j)\nonumber\\
=&0.\nonumber
\end{align}
\end{proof}
\end{thm}

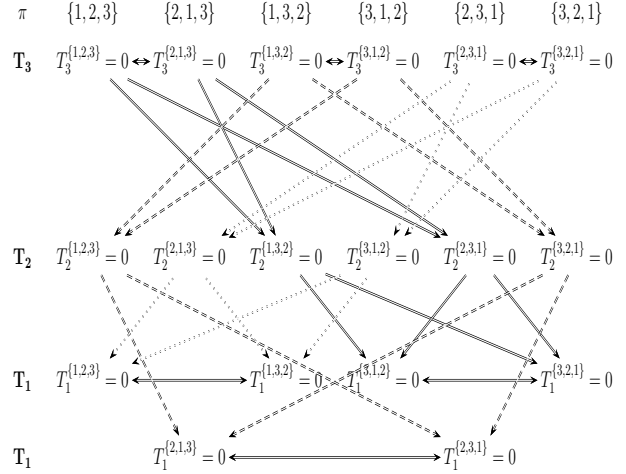
\begin{figure}[t]
  \noindent\resizebox{\linewidth}{.75\linewidth}{\input{exmpl4.1}}
  \caption{In the implication scheme above, the two-sided implications shown by $\Longleftrightarrow$ (e.g. $T^{\{1,2,3\}}_3=0~\Longleftrightarrow~T^{\{2,1,3\}}_3=0$, or $T^{\{1,2,3\}}_1=0~\Longleftrightarrow~T^{\{1,3,2\}}_1=0$) follow from Clause 1. of Theorem \ref{thm:genExp} without any extra assumptions.  The one-sided implications (e.g. $T^{\{1,2,3\}}_3=0~\Longrightarrow~T^{\{1,3,2\}}_2=0~\Longrightarrow~T^{\{3,1,2\}}_1=0$) follow from Clause 2 of the same theorem, and require specific posterior independence conditions to hold.  In particular, those implications that require the relation $\cind{V_3}{(V_1, V_2)}{\mc{D}}$ to hold are denoted by plain right arrows. On the other hand, those require $\cind{V_2}{(V_1,V_3)}{\mc{D}}$ are denoted by dashed right arrows.  Finally the implications that hold under $\cind{V_1}{(V_2,V_3)}{\mc{D}}$ are shown by dotted right arrows.}
  \label{fig:thm}
\end{figure}

\begin{exm} Let $\mc{V}=\left\{V_1,V_2,V_3\right\}$, ie. $I=\{1,2,3\}$. Suppose the $4$-term expansion of PPV is given by:
\begin{align}
  ~&{\sf Var}\left(Y_{n+1}\vert \mc{D},V_1,V_2,V_3\right)\nonumber\\
  =&T^{\{1,2,3\}}_0+T^{\{1,2,3\}}_3+T^{\{1,2,3\}}_2+T^{\{1,2,3\}}_1.
\end{align}
The implications of Theorem \ref{thm:gStCond} are displayed in Figure \ref{fig:thm} above. 
 \qed
\end{exm}

\subsection{More General Expansions of the PPV}

Let $\mc{V}_\mc{M}=\left\{V_1, V_2,\ldots, V_M\right\}\subseteq\mc{Z}$ be
an ordering of the set of variables in the model. Assume 
that $\mc{V}_{\mc{M}}$ is manifest and the rest of the 
variables in $\mc{V}\setminus\mc{V}_{\mc{M}}$ 
are latent. Consider an expansion of the PPV in terms of the elements of $\mc{V}_{\mc{M}}$.  
For any $u\le M\le K$ a $(u+1)$-term expansion of $Var\left[Y_{n+1}\vert \mc{D}\right]\left(\mc{V}_{\mc{M}}\right)$ can be achieved by first splitting the set $\mc{V}_{\mc{M}}$ into $u$ mutually exclusive and collectively exhaustive subsets, and then applying the LTV $u$-times to the subsets.  
The expressions for these expansions follow directly 
from \eqref{Conditional_Var_sum} apart from allowing each conditioning set of variables to have more than one
element.  

Clearly, for any given $u>1$, multiple expansions are 
possible, depending on the permutation of the subsets.  
Theorems \ref{thm:gStCond} and \ref{thm:genExp} will also hold mutatis mutandis.  
Our next result counts the number of 
expansions in the C-scope; recall the discussion 
after Example \ref{exm:hmm}.

\begin{proposition}
\label{Cscopecard1}
Fix $K$, $M$, $u$ as defined above, and suppose $u \leq M \leq K$. Then the
number of possible $(u+1)$-term PPV expansions of 
$Var\left[Y_{n+1}\vert \mc{D}\right]\left(\mc{V}_{\mc{M}}\right)$ is
\begin{align}
u!{{K}\choose{M}}S(M, u).
\label{Stirling}
\end{align}
Consequently, for fixed $K$, the total number of expansions of the PPV is given by
\begin{align}
\sum_{u=1}^K \sum_{M=u}^K u! {{K}\choose{M}} S(M, u).
\label{totalStirling}
\end{align}
Here, $S(M, u)$ is the Stirling number of the second kind with fixed $u\le M$.

\begin{proof}
By definition, $S(M, u)$ is the number of ways to form non-void, disjoint, and exhaustive collections of $u$ subsets from $M$ distinct objects. 

We start by observing that for a $(u+1)$-term expansion of the PPV, we must have $M$
disjoint non-void subsets of $(V_1, \ldots , V_K)$, i.e., not counting permutations, there are $S(M, u)$ possible choices.  Since we can permute these sets any way we want, we get a factor of $u!$.  Since we can do this for any choice of $M$ manifest variables out of $K$ variables, we get \eqref{Stirling}.  Summing over all the possible values of $M$and $u$ gives \eqref{totalStirling}.
\end{proof}

\end{proposition}

\begin{exm} For $K=2$, from \eqref{totalStirling} there are five possibilities that can be listed as follows.  For $K=2$, $M$ could be either $1$, or $2$.  For $M=1$ and $u=1$, there are two possibilities:  either $V_1$ alone or $V_2$ alone is manifest, that is, $V_2$ or, respectively, $V_1$ is latent. 
For $M=2$ and $u=1$, there is one possibility,
condition on $(V_1,V_2)$.  For $M=2$ and $u=2$ there are two possibilities:
condition on $V_1$ and then $V_2$ or condition on $V_2$ an then $V_1$. 
\end{exm}

\section{Term-wise Relative Size of Uncertainty}
\label{Draper}

In this section, we examine the relative sizes of the terms in PPV expansions. 
We look at a two-way random coefficient ANOVA model and then
revisit a computed example first treated in \cite{Draper:1995}.  Both
of these are three-term expansions.

We begin with an observation on asymptotic relative sizes of terms, 
borne out in the examples of Sec. \ref{sec:2}
It is easy to see that under standard 
regularity conditions (independent data with well-behaved conditional
densities and
$V_k$'s having well-behaved distributions)
the leading EVar term in two-term expansions is ${\cal{O}}_p(1)$, often simply
a constant given by the integral of a variance function.  Moreover, the 
second term is typically
going to be ${\cal{O}}_p(1/n)$ for reasons of posterior normality.  For 
three term expansions -- which actually
covers all cases simply by concatenating individual $V_k$'s on each side of
the variance operation --  we find that, under regularity conditions, the
leading term is ${\cal{O}}_p(1)$ and the last term is ${\cal{O}}_p(1/n)$.
However, the middle term is either ${\cal{O}}_p(1)$ if the inner variance
does not depend on the data, i.e., is only a function of $V_1$, or 
${\cal{O}}_p(1/n)$ if posterior normality holds for $V_2$ for each $V_1=v_1$.
That is, a PPV breaks down into two parts; the first part is the terms that
are ${\cal{O}}_p(1)$ and the second part is the terms that are ${\cal{O}}_p(1/n)$.

Next we derive expressions for the terms in two-way random coefficient model
so that we can graph them.  After that we re-analyse an example from
\cite{Draper:1995}.

\begin{figure*}[h!]
\subfigure[Value of the terms. \label{fig:termT}]{
  \resizebox{.75\columnwidth}{!}{\includegraphics{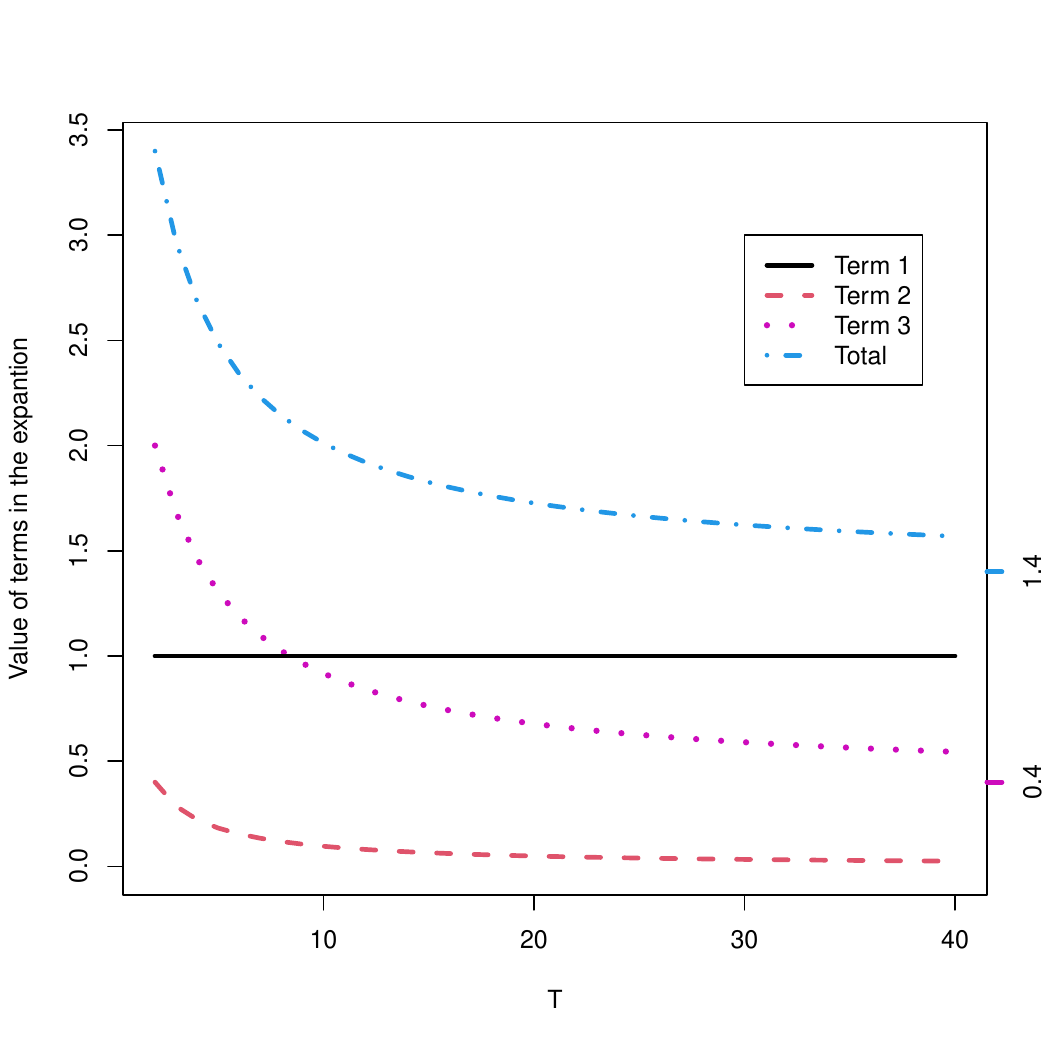}}
  }\qquad\subfigure[Proportion of the terms to the total. \label{fig:propT}]{
  \resizebox{.75\columnwidth}{!}{\includegraphics{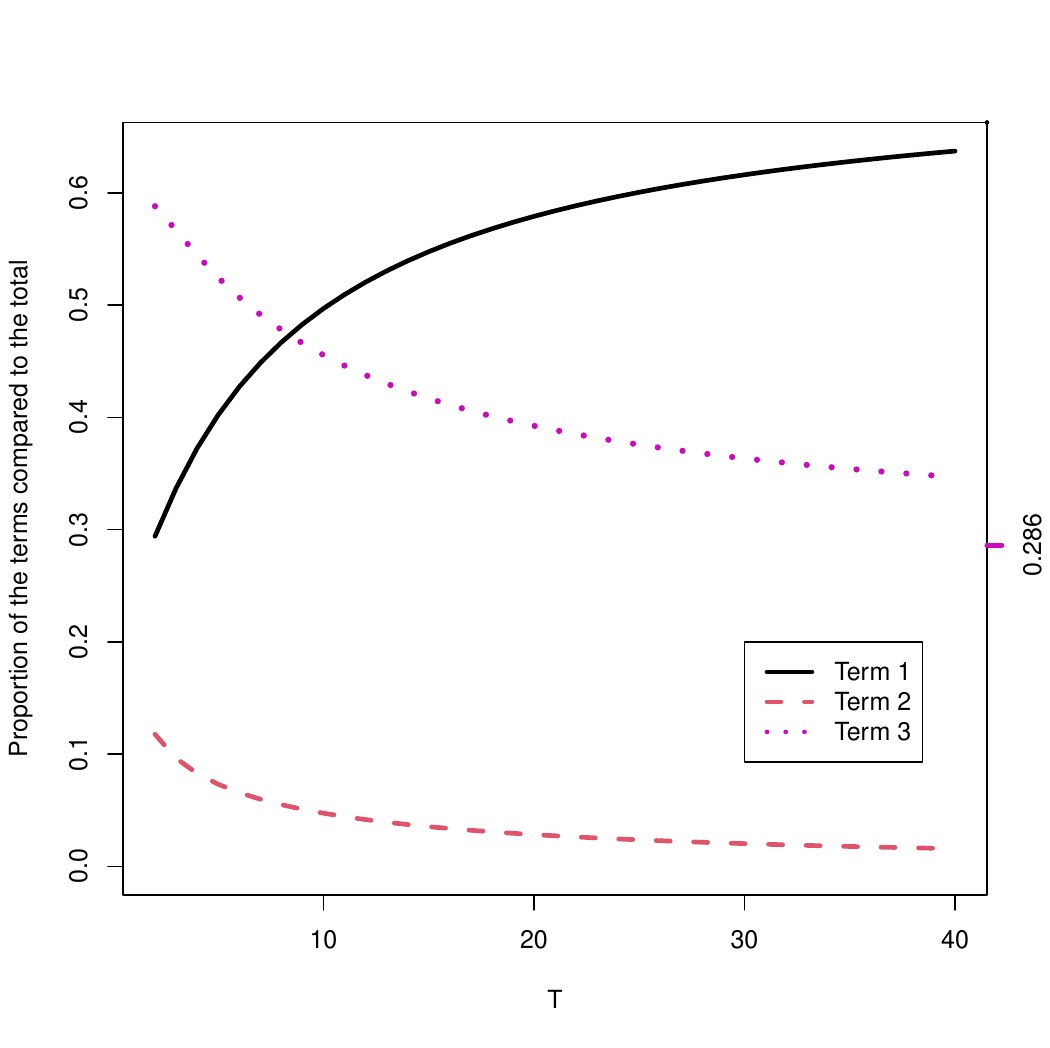}}
  }
\subfigure[Value of the terms. \label{fig:termB}]{
  \resizebox{.75\columnwidth}{!}{\includegraphics{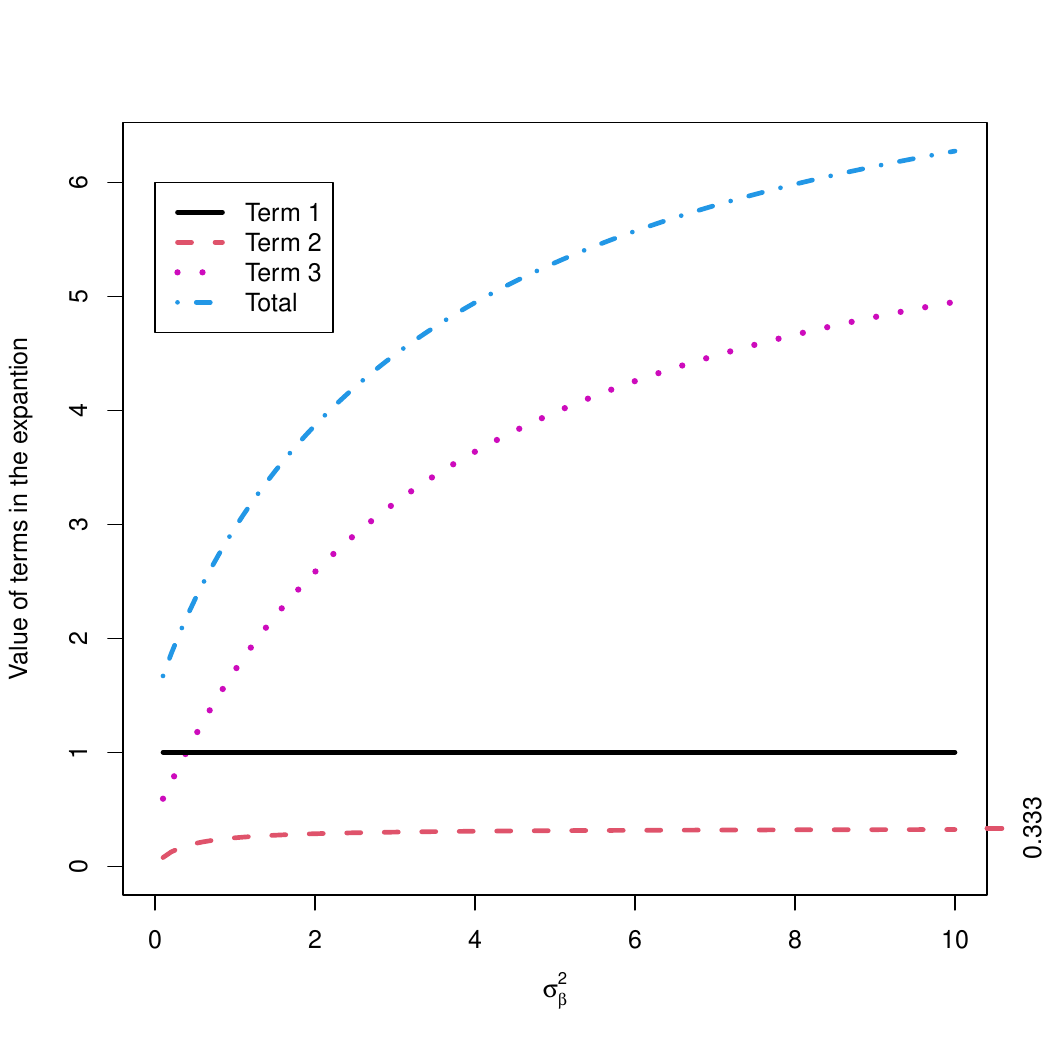}}
  }\qquad\subfigure[Proportion of the terms to the total. \label{fig:propB}]{
  \resizebox{.75\columnwidth}{!}{\includegraphics{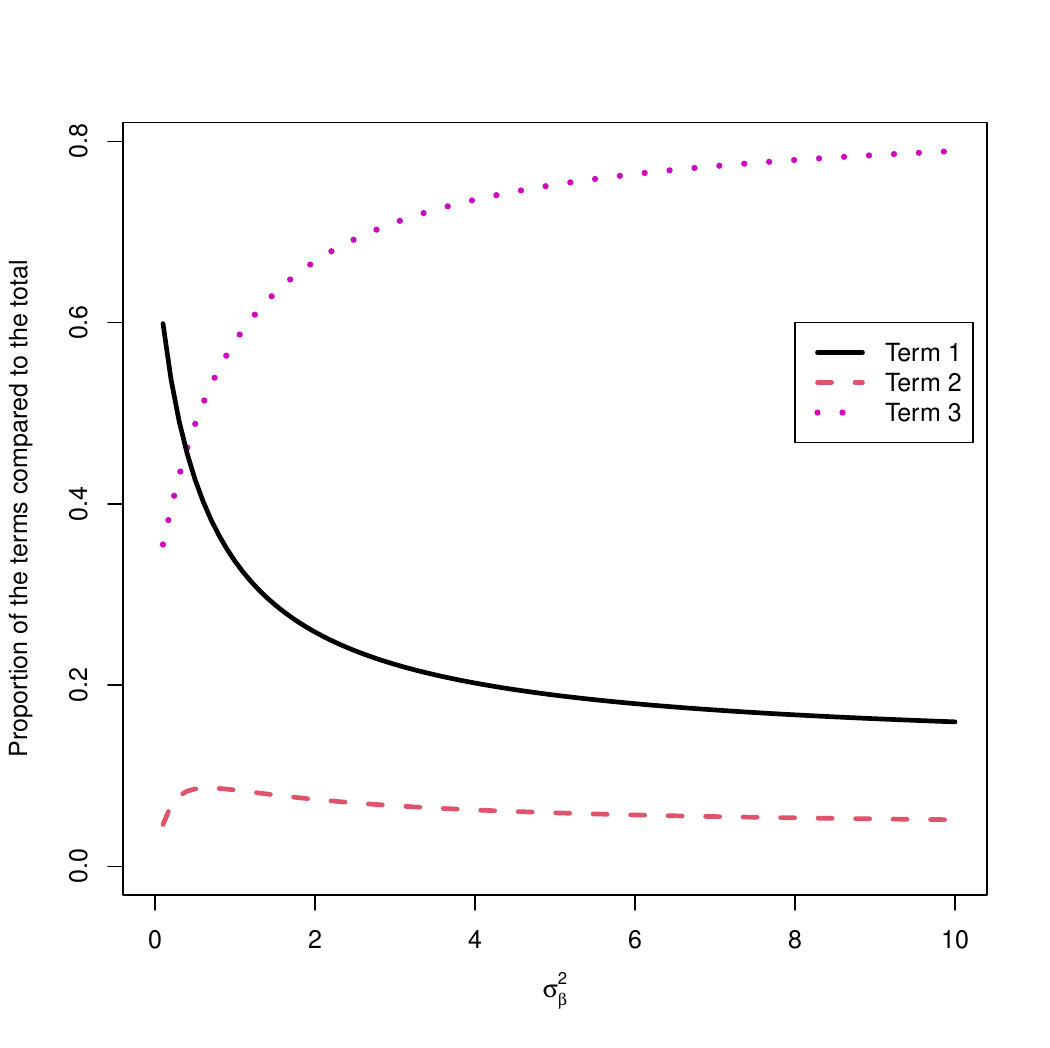}}
}
\caption{The value and the proportion of the three terms in \ref{2wayANOVAvardecomp0}.  Figures \ref{fig:termT} and \ref{fig:propT} are for $\sigma^2_{\tau}=\sigma^2_{\beta}=2$, $\sigma^2_{\epsilon}=1$, $B=2$, and varying $T$.  For Figures \ref{fig:termB} and \ref{fig:propB}, we have set $B=2$, $T=3$, $\sigma^2_{\tau}=5$, $\sigma^2_{\epsilon}=1$, and let $\sigma^2_{\beta}$ vary.}
\label{2wayANOVAgraphs}
\end{figure*}

\subsection{Two-Way Random Coefficient Model}
\label{2wayANOVAeg}

Consider the Bayesian model defined by
\begin{eqnarray}
Y_{ij} = \tau_i + \beta_j + \epsilon_{ij},
\label{2wayANOVAmodel}
\end{eqnarray}
where $i=1, \ldots , T$, $j=1, \ldots , B$ and we set
\begin{eqnarray}
\tau_i &\sim& N(\tau_0, \sigma^2_\tau)
\nonumber \\
\beta_j &\sim& N(\beta_0, \sigma^2_\beta)
\nonumber \\
\epsilon_{ij} &\sim& N(0, \sigma^2_\epsilon)
\label{distributions}
\end{eqnarray}
with
$$
\tau_i  \upmodels \tau_j ,\quad \beta_i \upmodels \beta_j 
$$
for $i\neq j$ and for all $i$, $j$
$$
\tau_i \upmodels \beta_j , \quad  \tau_i, \beta_j \upmodels \epsilon_{ij};
$$
the notation $U\upmodels V$ means that $U$ is marginally independent of $V$.
Our data is of the form ${\bf y}=\left(y_{11},\ldots,y_{TB}\right)$, where 
$n=T\times B$ and we want to predict a future $Y_{ij}$ for 
some prespecified $i$ and $j$ with
$1\le i\le T$ and $1\le j\le B$.

Because of \eqref{2wayANOVAmodel}, the future $Y_{ij}$ will depend on the entire matrix ${\bf y}$.  
Thus, the predictor automatically borrows information from ``similar'' observations.   So, PI's will provide a quantification 
of the uncertainty for the $ij$-th cell averaging over the parameters in the model.  Models such as 
\eqref{2wayANOVAmodel} arise naturally in small-area estimation, poverty mapping, missing plot analysis in 
agronomy, and Bayesian imputation for missing values in contingency tables amongst other settings.

The three-term expansion of the PPV, with $V_1=\beta$ and $V_2=\tau$ is
\begin{align}
 ~&{\sf Var}(Y_{ij} \vert {\bf y})\nonumber\\
=& {\sf E}_{\tau \vert {\bf y}} {\sf E}_{\beta \vert  {\bf y}, \tau } {\sf Var}(Y_{ij} \vert {\bf y}, \beta, \tau)
\nonumber \\
& + {\sf E}_{\tau \vert {\bf y}}  {\sf Var}_{\beta \vert  {\bf y}, \tau } {\sf E}(Y_{ij} \vert {\bf y}, \beta, \tau)
\nonumber \\
& + {\sf Var}_{\tau \vert {\bf y}}  {\sf E}_{\beta \vert {\bf y}, \tau } {\sf E}(Y_{ij} \vert {\bf y}, \beta, \tau) .
\label{2wayANOVAvardecomp0}
\end{align}

From the detailed derivations in Appendix \ref{calcs2wayANOVA} we get


\[
  {\sf Term } ~ 1 = \sigma^2_\epsilon, {\sf Term }  ~ 2 =  \left( \frac{T}{\sigma_\epsilon^2} + \frac{1}{\sigma^2_\beta}  \right)^{-1}, \quad \mbox{and}
  \]
  
\smallskip
\begin{widetext}
\begin{align}
{\sf Term} ~3=& \frac{1}{
a\cdot\left( \frac{T}{\sigma^2_\epsilon} + \frac{1}{\sigma^2_\beta} \right)^2
} 
 \left\{ 
  \left( \frac{T+1}{\sigma^2_\epsilon} + \frac{1}{\sigma^2_\beta} \right)^2\left( 1 - \frac{Bb}{a+bBT}\right)+
  \frac{(T-1)}{\sigma^4_\epsilon}\left( 1 - \frac{Bb}{a + bBT} \right)\right.
  \nonumber \\
 &\quad - \left. 
 \frac{2}{\sigma^2_\epsilon} \left( \frac{T+1}{\sigma^2_\epsilon} + \frac{1}{\sigma^2_\beta} \right) (T-1) 
 \frac{Bb}{(a+ BbT)}
  \right\},
\nonumber
\end{align}
\end{widetext}

where $a = \frac{B}{\sigma^2_\epsilon} + \frac{1}{\sigma^2_\tau}$ and
$b = -  \left\{ \sigma^4_\epsilon \left( \frac{T}{\sigma^2_\epsilon} + \frac{1}{\sigma^2_\beta}\right) \right\}^{-1}$.  
The three terms are symmetric in the sense that if we interchange the order
of conditioning in \eqref{2wayANOVAvardecomp0} we will interchange the
the $\tau$'s and $\beta$'s in the terms.

The key issue is the relative sizes of these three terms.  Term 1 is easy to visualize because it's a constant.
Term 2 is next easiest because it does not depend on $B$ or $\sigma^2_\tau$.  It is seen that
Term 2 is smaller than both $\sigma^2_{\epsilon}$ and $\sigma^2_{\beta}$.

Indeed, for fixed $\sigma^2_{\epsilon}$ and $\sigma^2_{\beta}$, Term 2 is asymptotically $\mathcal{O}(1/T)$
and in the limit of $\sigma^2_{\beta}\rightarrow\infty$, Term 2 converges to 
$\sigma^2_{\epsilon}/T>0$ -- which is the variance of $\bar{y}_{\cdot j}$.

The behaviour of Term 3 is more complicated.  It depends on $B$ through $a$ and $T$ through $b$.
By comparing the numerator and the denominator, we see that, asymptotically for large $T$ holding the other constants fixed,
Term 3 is $\mathcal{O}(1)$.  
In the limit as $T\rightarrow\infty$, Term 3 convegres to 
$a^{-1}=(\sigma^2_{\epsilon}\sigma^2_{\tau})/(B\sigma^2_{\tau}+\sigma^2_{\epsilon})$.  
Thus, for large values of $T$, Term 3 is independent of $\sigma^2_{\beta}$ and
smaller than both $\sigma^2_{\epsilon}$ and $\sigma^2_{\tau}$.
Analogous reasoning applied to Term 3 shows that it is $\mathcal{O}(1/B)$ as $B$ increases,
when the other quantities are held constant.

It is now clear that unless both $B$ and $T$ diverge to infinity, Terms $2$ and $3$ won't simultaneously go to zero.
That is, if only one of $B$ or $T$ goes to infinity, at least one of Term 2 or 3 may be a significant fraction of
the PPV. 

When $T\rightarrow\infty$, Term $2$ $\rightarrow 0$, and the limiting PPV is
  \[
  \sigma^2_{\epsilon}+\frac{1}{a}=\sigma^2_{\epsilon}\frac{(B+1)\sigma^2_{\tau}+\sigma^2_{\epsilon}}{B\sigma^2_{\tau}+\sigma^2_{\epsilon}}.
\]
That is, under this asymptotic regime, relative size of Term $1$ to the total is $\{B\sigma^2_{\tau}+\sigma^2_{\epsilon}\}/\{(B+1)\sigma^2_{\tau}+\sigma^2_{\epsilon}\}$, and that of Term $3$ is $\sigma^2_{\tau}/\{(B+1)\sigma^2_{\tau}+\sigma^2_{\epsilon}\}$.

A similar phenomenon can be observed when $B\rightarrow\infty$.  
In this case Term $3$ $\rightarrow 0$, and the limiting PPV is
  \[
  \sigma^2_{\epsilon}\frac{(T+1)\sigma^2_{\beta}+\sigma^2_{\epsilon}}{T\sigma^2_{\beta}+\sigma^2_{\epsilon}},
  \]
which implies that, the relative size of term $1$ to the total is $\{T\sigma^2_{\beta}+\sigma^2_{\epsilon}\}/\{(T+1)\sigma^2_{\beta}+\sigma^2_{\epsilon}\}$ and that of Term $2$ to the limiting PPV is $\sigma^2_{\beta}/\{(T+1)\sigma^2_{\beta}+\sigma^2_{\epsilon}\}$,
i.e., neither Term 1 nor Term 2 can be neglected in general.

An illustration of the behaviour of the three terms is in Figure \ref{2wayANOVAgraphs}.  
In Panel \ref{fig:termT} we keep $\sigma^2_{\tau}=\sigma^2_{\beta}=2$, $\sigma^2_{\epsilon}=1$, $B=2$, and vary $T$.  For values of $T$ larger than $8$, Term 1 dominates the other two terms.  Term 2 will reduce to zero as $T\rightarrow\infty$.  Term 3 reduces 
to $a^{-1}=0.4$.  The PPV itself converges to $1.4$. In Panel \ref{fig:propT}, it can be seen that, 
around $T=20$,  Term 2 can be omitted at a threshold of about $5\%$.
In Panel \ref{fig:termB}, ($B=2$, $T=3$, $\sigma^2_{\tau}=5$, $\sigma^2_{\epsilon}=1$) Term 3 dominates, and Term 2 
converges to $0.33$, as $\sigma^2_{\beta}\rightarrow\infty$.  In Panel \ref{fig:propB} as $\sigma^2_{\beta}\rightarrow\infty$, 
all terms have limits in $(0,1)$.

In the above we have assumed we have only one replication of ${\bf Y}=\left(Y_{11},\ldots,Y_{T,B}\right) = {\bf y}$ to predict
the next outcomes.   More generally,
our reasoning can be extended to $n$ replications.  The corresponding expressions for
the terms in \eqref{2wayANOVAvardecomp0} will be more complex.

\subsection{Application to the Challenger disaster data}
\label{challenger}

For the Challenger Space Shuttle disaster, it is widely believed that making
the decision to launch the space shuttle at an ambient temperature and 
pressure at which various components had not been tested ended up being 
catastrophic.  
However, the disaster could have been avoided had a proper 
uncertainty analysis been done.  
Statistically, the error of the decision makers was to choose a single model from a
model list rather than incorporating all sources of predictive 
uncertainty into their analysis.

We apply the PPV expansion to the Challenger disaster data, where the probability 
of O-ring failure on a space shuttle is predicted for a given temperature and 
pressure.  This data set was analyzed by \cite{Draper:1995}, who pointed 
out the need for evaluating the entire uncertainty associated with what he called
`structural' choices.  Below, we see that the relative sizes of the terms 
in the three-term PPV expansion depend heavily on the $V_k$'s
used in the analysis.

The data consists of $22$ observations (one observation out of $23$ is 
generally considered an outlier) of the number of damaged O-rings.  
The shuttle had six primary O-rings. The maximum number of recorded failures 
was two. Together with the number of failures, the temperature $t$ of the 
joint, and the field leak-check pressure $s$ were recorded. More 
details of the data can be found in \cite{Dalal:etal:1989}.  

It is widely believed that the O-ring resiliency was directly related to its temperature. Experimental data showed: ``At $100^{\circ}$F, the O-ring maintained contact. $\ldots$ At $50^{\circ}$F, the O-ring did not reestablish contact" (see \cite{Dalal:etal:1989}).  The ambient temperature on the day of launch was $31^{\circ}$F, and the recorded field leak-check pressure was $200$ psi (see \cite[Table 1]{Dalal:etal:1989}).  The general goal in the literature has been to predict the probability of failure of one or more O-rings on the day of launch. It is generally predicted that the probability of failure was quite high, close to one.    

In order to analyse the PPV of the predicted failure probability,  
we assume that the number of damaged O-rings follows a $Binomial(6, p_{t,s})$ 
distribution
where $p_{t,s}$ is the probability of failure of one O-ring at temperature $t$ and pressure $s$.  The probability $p_{t,s}$ is connected to the explanatory 
variables via one of three link functions, namely, logit, $c\log\log$, and probit.  

In their analysis, \cite{Dalal:etal:1989} used several models with temperature and pressure as covariates.  They argued that the temperature is the most important variable; however, a weak effect of leak-check pressure could not be categorically rejected.  This may imply that some latent variables weakly correlated to pressure are missing from the model. However, for prediction purposes, we don't 
necessarily need the true model.  Here, we include temperature $t$ and pressure 
$s$ in our model and, following  \cite{Draper:1995}, we also use $t^2$.
We expand our predictive variance in terms of two conditioning
variables, namely the link function $V_1= \{$logit, c$\log\log$, probit$\}$
and $V_2=$ the set of all $15$ models that can be formed by choosing variables from 
$\{1, t , t^2, s\}$.  Hence, we use a three-term expansion of the PPV.

By contrast, by making some simplifying assumptions leading to fewer
models than we use here, \cite{Draper:1995} 
provided a table of 
posterior quantities for the `structural' choices, and a two-term 
PPV expansion for `within-structure' and `between-structure' variances as
\begin{align}
  ~&{\sf Var}(p_{t=31}\vert  {\cal{D}}_{23})= {\sf Var}_{within} + {\sf Var}_{between}\nonumber\\
  =& 0.0338 + 0.0135 = 0.0473.
\label{twotermDraper}
\end{align}
Because
$.0135/.0473 \approx 28.5\%$, Draper concluded that the uncertainty 
represented by the second term
in \eqref{twotermDraper} could not be neglected, thereby giving a much larger
PI for $p_{t=31}$ including unacceptably large probabilities of failure.

\begin{figure}[t]
    \resizebox{2.5in}{!}{\includegraphics{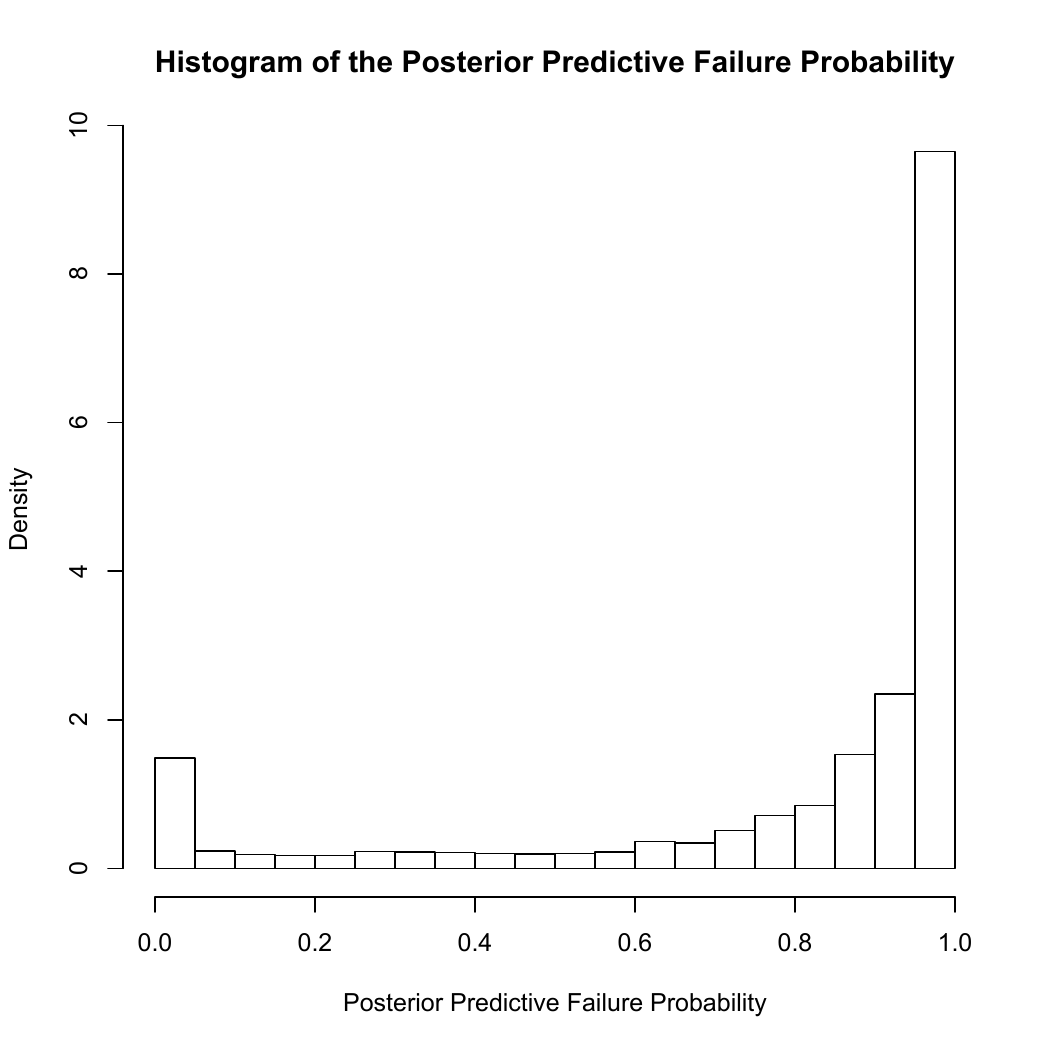}}
    \caption{Histogram of the posterior predicted failure probabilities $p_{t=31,s=200}$ for the Challenger space shuttle on the day of its launch. }
    \label{fig:predhist}
\end{figure}

Here, we calculated  posterior predicted failure probabilities $p_{t=31,s=200}$ 
given the data by
sampling using a reversible jump Markov Chain Monte Carlo \citep{green1995}, 
with regression coefficients following a (diffuse) $N\left(0,100\right)$ prior, 
and assuming uniform priors on $V_1$ and $V_2$.  The results presented are 
based on $10$ million draws with a burn-in of $5$ million.  

The histogram of the predicted posterior probabilties is in Figure
\ref{fig:predhist}.  The histogram clearly shows that the posterior prediction failure probabilities are large and are close to one.  That is, the 
failure of the O-rings should have been anticipated.  From the histogram,
we also see that the posterior predictive variance is large. 

We can 
quantify the sources of uncertainty using three-term decompositions 
with respect to the models and the link functions.  
The three-term expansion corresponding to \eqref{LTV2termgen} is
\begin{align}
\label{var_challenger_3term}
~&{\sf Var}(p_{t=31,s=200}\vert  {\cal{D}}_{23})\nonumber\\
=&  {\sf E}_{V_1}{\sf E}_{V_2}{\sf Var}(p_{t=31,s=200} \vert  {\cal{D}}_{23}, V_1, V_2)\nonumber\\
&+ {\sf E}_{V_1}{\sf Var}_{V_2}E(p_{t=31,s=200} \vert  {\cal{D}}_{23}, V_1, V_2)\nonumber\\
& + {\sf Var}_{V_1}{\sf E}(p_{t=31,s=200} \vert  {\cal{D}}_{23}, V_1)\nonumber \\
=& 0.07212 + 0.01682 + 0.00636=0.09531.
\end{align}
Furthermore, by permuting the position of $V_1$ and $V_2$, the three-term expansion corresponding to \eqref{LTV2termgen2} becomes:
\begin{align}
\label{var_challenger_3termb}
~&{\sf Var}(p_{t=31,s=200}\vert  {\cal{D}}_{23})\nonumber\\
=& {\sf E}_{V_2}{\sf E}_{V_1}{\sf Var}(p_{t=31,s=200} \vert  {\cal{D}}_{23}, V_1, V_2)\nonumber\\
&+ {\sf E}_{V_2}{\sf Var}_{V_1}{\sf E}(p_{t=31,s=200} \vert  {\cal{D}}_{23}, V_1, V_2)\nonumber\\
& + {\sf Var}_{V_2}{\sf E}(p_{t=31,s=200} \vert  {\cal{D}}_{23}, V_2)\nonumber\\
=& 0.07212 + 0.005 + 0.01819=0.09531.
\end{align}

In both expansions, the first term dominates. However, the second largest
terms -- that take a variance over $V_2$  -- clearly dominate the terms 
involving variances over $V_1$.  This is reasonable because the uncertainty due to
the model i.e., $V_2$, is intuitively higher than the uncertainty due
to the link function, i.e., $V_1$.  Indeed, both expansions suggest we can drop
$V_1$ from our analysis leaving a two term PPV.   Moreover, if 
the hypothesis of Lemma \ref{condind} holds, at least approximately, then we
get Theorem \ref{thm:genExp} at least approximately.  So, we surmise
term 3 in \eqref{var_challenger_3term} is near zero because
term 2 in \eqref{var_challenger_3termb} is.

This example shows that examining the terms in an expansion of the PPV
can lead us to identify aspects of a model that just don't matter.
Indeed, here, we could simply choose one of the link functions and
use it, ignoring the other `levels' of $V_1$.  This is analogous to
testing whether a factor in a $K$-way ANOVA can be left out.
The obvious way to formalize this is by choosing a user specified level,
say $\alpha$, and omitting levels in a HM that contribute less than $\alpha\%$
of the PPV.

\section{Discussion}
\label{discuss}

The main contribution of this paper is to provide an
expansion of the posterior predictive
variance (PPV) for hierarchical models (HM's) and models that, more generally, satisfy conditional independence assumptions that we have
called structural.   An immediate benefit from this is that
we have a conservation law over expansions for the PPV.   This is important for two
reasons.  First, the PPV controls the width of prediction intervals 
so we want to know what aspects of variance are contributing 
most to it.    Second,
we want to identify what levels of a model can be 
collapsed to a single value.  This is
analogous to testing for whether a factor can be 
dropped in a frequentist multi-way ANOVA.

Our expansions start with a fixed hierarchical model and 
hence a fixed PPV that can be expressed in multiple 
expansions depending on the how the law of
total variation (LTV) is iterated.
The various expansions depend on the
ordering of the conditioning variables from the levels of the model.
We focus on what we call the $C$-scope of a model -- the
collection of expansions of the posterior predictive variance that arise from 
using the LTV
only on terms in which an expectation of a 
variance appears. In Sections \ref{sec:quant} and \ref{decomposition}
we give an extensive discussion of when a term being zero 
can be used to imply terms in other expansions
must be zero, too.  In Proposition
\ref{Cscopecard1} we give an explicit expression for the cardinality of the $C$-scope.

The main methodological implication of our work is that we can 
more readily use HM's where we might have used Bayesian 
nonparametrics.   Indeed, we can represent any feature of
a statistical model with a prior.  These features may or may not
have any physical analog: in Subsec. \ref{challenger} we use 
variable selection as a level in an HM and this is part of 
physical modeling.  We also take the link function in a 
GLM as a feature and this is not in general seen as an aspect of
modeling.  Elsewhere, see \cite{Dustin:Ghosh:Clarke:2025}, we used selection of a shrinkage method as a feature of modeling 
and this does not really have a physical analog; 
indeed it corresponds to selecting a prior.

One effect of using variance is that the metric properties of the 
model list become important as well as its probabilistic 
properties.   Thus, as a matter of model list design we 
want to choose a HM so that its PPV is not too small 
relative to the data so that using multiple 
expansions to prune out levels will be effective.   
We also want the PPV to be not too large or the
pruning will be obvious without analysis and could lead to bias.  We want the
elements of the model list to be close enough to each other
that they are plausible but far enough apart from each other that they are distinguishable with the data we have.  Indeed, we 
may want to construct a HM so that the higher the level the
less it is thought to matter and then order the uses of the 
LTV so that we start by conditioning
on the level we think will be easiest to eliminate. 

It is a sort of folk-theorem that the higher the 
level the less important is and
our procedure can assess this, see Theorem \ref{thm:var}. 
Even though we can
construct examples of arbitrarily many levels in which the top
level does matter, the intuition holds and extensions of our 
work may be able to provide a formal
way to decide if upper levels in a HM should be retained.

One drawback of our procedure is that it we do not have 
a fully formal way to assess the relative contributions
of terms in the expansion.   We have relied on essentially a 
user specified threshold for
whether a term is large enough to retain.  
This is so because in general we do not have a likelihood
for these terms and therefore cannot do Bayes testing directly,
although frequentist tests are possible, see \cite{Dustin:Ghosh:Clarke:2025}.   
On the other hand,
there are ways around this e.g.,  pseudo-Bayes posteriors in 
which a likelihood is formed from
an empirical risk.  
We have not investigated this
possibility, but it is promising as it is in the 
spirit of the mathematical modeling we
advocate here, namely being willing to use 
mathematical quantities without
physical motivation as a way to produce predictive analyses.

We conclude with two final observations.
First,  the interplay between variance and conditioning is not well understood because 
variance is non-linear.   Our methods provide a context where this interplay can be explored.
Moreover, as we can see, each level of modeling can be
treated as a level in an uncertainty quantification.  So, given an ordering on the components in a
HM, our analysis here is a step towards
assigning an uncertainty to each component of the model.

Second, the treatment we have given for variance can, in 
principle, be extended
to other nonlinear operations, though it looks hard.  
For instance, \cite{Brillinger:1969} 
gives a way to calculate 
cumulants of a distribution that can be a posterior quantity.
He gives a formula similar to  
\ref{Conditional_Var_sum} and gives examples
using this result for sums of variables and mixture distributions. 
In addition, we could have used the Shannon mutual information in 
place of the variance and invoked its chain rule.  We have not 
chosen these because the first seems
quite hard and the second is not as readily applicable to 
prediction intervals.


\begin{appendix}

\section{Calculations for the Three Term Normal}
\label{calcs3termnormal}

Our task is to derive an expression for 
${\sf Var}(Y_{n+1} \vert y^n)$ directly.
Recall \eqref{LTV2termgen} and that
we have two parameters $\mu$ and $\lambda$
as well as three hyperparameters $\kappa_0$, $\alpha_0$, and $\beta_0$.  For simplicity, write
$\gamma = \lambda^2$.   The conditional density of $y^n$ given $\mu$ and $\lambda^2$ is
\begin{align}
~&p(y^n \vert \mu, \gamma)\nonumber\\ 
=& 
\gamma^{n/2} e^{-\gamma/2 \sum_{i=1}^n (y_i - \mu)^2}
\sqrt{\gamma \kappa_0}  e^{- \kappa_0 \gamma/2 \sum_{i=1}^n (\mu - \mu_0)^2}\nonumber\\
&\times\frac{\beta_0^{\alpha_0}}{\Gamma(\alpha_0)} \gamma^{\alpha_0 -1} e^{- \beta_0 \gamma} .
\label{likelihood}
\end{align}
We have that
\begin{align}
~&p(y_{n+1} \vert y^n)= \int \int p(y_{n+1} \vert y^n, \mu, \gamma)p(\mu, \gamma \vert y^n) {\rm d} \mu {\rm d} \gamma
\nonumber \\
=&    \int \int p(y_{n+1} \vert y^n, \mu, \gamma)p(\mu \vert  y^n, \gamma )  p(\gamma \vert y^n){\rm d} \mu {\rm d} \gamma.
\label{actualpreddty}
\end{align}
We want to identify the three densities in the integrand.  We know the first.

For the second, with some foresight, let
\begin{eqnarray}
\begin{aligned}
\mu_n &=  \frac{n\bar{y} + \kappa_0\mu_0}{\kappa_n}
\\
\gamma_n&= \gamma(n+\kappa_0) = \gamma \kappa_n
\\
\kappa_n &= n + \kappa_0
\\
\alpha_n &= \alpha_0 + (n/2)
\\
\beta_n &= \beta_0 + \frac{1}{2} \left( \sum_{i=1}^n y_i^2 - \frac{\kappa_0 \mu_0 - n \bar{y}}{\kappa_0+ n} \right).
\end{aligned}
\label{defsNormal}
\end{eqnarray}

\noindent{\it Step 1:}
We begin by seeing that
\[
p(\mu \vert y^n, \gamma) \sim N(\mu_n, 1/\gamma_n).
\]
The squared terms in the exponent in \eqref{likelihood} are
\begin{align}
-&\frac{\gamma}{2} \sum_{i=1}^n(y_i - \mu)^2 - \frac{\kappa_0 \gamma}{2} (\mu - \mu_0)^2
\nonumber \\
=&
- \frac{\gamma}{2} \bigg[ \sum_{i=1}^n y_i^2 + n \mu^2 - 2 n \bar{y}\mu + \kappa_0 \mu^2\nonumber\\
  &\hspace{1in}+ \kappa_0 \mu_0^2 - 2 \kappa_0\mu \mu_0\bigg]
\nonumber \\
=&
- \frac{\gamma}{2} \bigg[ \mu^2(n + \kappa_0) - 2 \mu (n\bar{y} + \kappa_0 \mu_0)\nonumber\\
 &\hspace{1in} + \sum_{i=1}^n y_i^2 + \kappa_0\mu_0^2 \bigg]
\label{completesquare1}
\end{align}
Completing the square in $\mu$ means \eqref{completesquare1} becomes

\begin{align}
-&\frac{\gamma}{2} 
\left[ \mu^2 (n + \kappa_0) - 2 \mu \sqrt{n+\kappa_0} \frac{(n\bar{y} + \kappa_0 \mu_0)}{\sqrt{n+\kappa_0}}\right.\nonumber\\
 &\hspace{1in}\left. + \frac{(n \bar{y}+ \kappa_0 \mu_0)^2}{n+ \kappa_0} \right]
\nonumber \\
& \quad -
 \frac{\gamma}{2} 
\left[ \sum_{i=1}^n y_i^2 + \kappa_0\mu_0^2 - \frac{(n \bar{y} + \kappa_0 \mu_0)^2}{n+ \kappa_0} \right]
\nonumber \\
=&
- \frac{\gamma(n+\kappa_0)}{2}
\left[ \mu - \frac{n\bar{y} + \kappa_0 \mu_0}{n +\kappa_0} \right]^2\nonumber\\
&-\frac{\gamma}{2}
\left[ \sum_{i=1}^n y_i^2 + \kappa_0 \mu_0^2 - \frac{(n\bar{y} + \kappa_0 \mu_0)^2}{n +\kappa_0}  \right].
\label{completesquare2}
\end{align}


Note that the `extra' $\sqrt{\gamma}$ in \eqref{likelihood} is absorbed into the normal density.
This completes Step 1.

\noindent{\it Step 2:}  Next, we see that $p( \gamma \vert y^n) \sim {\sf Gamma}(\alpha_n , \beta_n)$.
By exponentiating the second term in \eqref{completesquare2} and multiplying it by the `active' factors 
in \eqref{likelihood} we get that the rest of  the likelihood is proportional to
\begin{equation}
\gamma^{\alpha_0 + (n/2) - 1} \gamma^{n/2} e^{-\gamma (\beta_0+ (1/2)\left[ \sum_{i=1}^n y_i^2 +\kappa_0 \mu_0^2
- \kappa_n \mu_n^2 \right] }.
\label{cdlpstr1}
\end{equation}
Upon normalization this gives Step 2.

We comment that in principle, we now have the right hand side of \eqref{actualpreddty}.  However,
finding ${\sf Var}(Y_{n+1} \vert y^n)$ directly is a lot of work (probably involving $t$-distributions).  
So, we use a two term expansion.  For this we derive the following.

\noindent{\it Step 3:}   Obtain the conditional posterior
\begin{eqnarray}
p(y_{n+1} \vert y^n, \gamma) \sim N\left(\mu_n, \frac{\kappa_n + 1}{\kappa_n \gamma}\right).
\label{cdlpstr2}
\end{eqnarray}
To do this, first note
\begin{eqnarray}
p(y_{n+1} , \gamma \vert y^n) = p(y_{n+1} \vert y^n, \gamma) p(\gamma \vert y^n).
\label{factor}
\end{eqnarray}
Since we have $p(\gamma \vert y^n)$ it is enough to find the right hand side of \eqref{factor}.  To do this
recall that by definition
\begin{align}
~&p( y _{n+1}, \gamma \vert  y^n)\nonumber\\   
=&p(\gamma \vert y^n) \int p(y_{n+1} \vert \mu, \gamma) p(\mu  \vert y^n, \gamma)  {\rm d} \mu 
\label{preddty}
\end{align}


The integrand in \eqref{preddty} (in $\mu$) is
\begin{eqnarray}
&& \propto \sqrt{\gamma} e^{-(\gamma/2)(y_{n+1} - \mu)^2 } \times  \sqrt{\gamma_n} e^{-(\gamma/2) (\mu - \mu_n)^2}
\nonumber\\
&& = \gamma \sqrt{\kappa_n} e^{-(\gamma/2) \left[  (y_{n+1} - \mu )^2 + \kappa_n (\mu - \mu_n)^2 \right]} .
\label{sumsquares}
\end{eqnarray}
By some notational gymnastics, completing the square in \eqref{sumsquares} gives that 
\begin{align}
~& (y_{n+1} - \mu )^2 + \kappa_n (\mu - \mu_n)^2 
\nonumber \\
~& =
(1 + \kappa_n) \left( \mu  - \frac{y_{n+1} + \kappa_n \mu_n}{1 +\kappa_n} \right)^2
+ y_{n+1}^2\nonumber\\
&+\quad \kappa_n \mu_n^2 -  \frac{(y_{n+1} + \kappa_n \mu_n)^2}{1 +\kappa_n}  .
\label{compdsquare}
\end{align}
Using  \eqref{compdsquare} in \eqref{sumsquares}
gives that the integrand in \eqref{preddty}  is
\begin{align}
  \propto& \frac{ \sqrt{\kappa_n} }{\sqrt{1+\kappa_n}}  \sqrt{\gamma (1 + \kappa_n) }\nonumber\\
  &\times e^{-(\gamma/2) 
\left[(1 + \kappa_n) \left( \mu  - \frac{y_{n+1} + \kappa_n \mu_n}{1 +\kappa_n} \right)^2 \right]}
\nonumber \\
&\times \sqrt{\gamma} e^{- (\gamma/2) \left[  y_{n+1}^2 + \kappa_n \mu_n^2 -  \frac{(y_{n+1} + \kappa_n \mu_n)^2}{1 +\kappa_n} \right]  } .
\label{integrand}
\end{align}
The first factor can be integrated over $\mu$ and the exponent in the second factor is
\begin{align}
& y_{n+1}^2 + \kappa_n \mu_n^2 -  \frac{(y_{n+1} + \kappa_n \mu_n)^2}{1 +\kappa_n} 
\nonumber \\
& \quad = \frac{1}{1 + \kappa_n} \left[ \kappa_n \left(y_{n+1}^2 +  \mu_n^2 - 2 y_{n+1}\mu_n. \right)  \right]
\nonumber \\
& \quad = \frac{\kappa_n}{1 + \kappa_n} \left( y_{n+1} - \mu_n \right)^2 .
\label{2ndfactor}
\end{align}
Now we see that the integral in \eqref{preddty} gives \eqref{cdlpstr2}, completing Step 3.

To complete the derivation of the PPV, write
\begin{align}
  ~&{\sf Var}(Y_{n+1} \vert y^n)\nonumber\\
  =& {\sf E}\left[ {\sf Var} (Y_{n+1} \vert y^n, \gamma) \right]
+ {\sf Var} \left[ {\sf E}(Y_{n+1} \vert  y^n, \gamma) \right]
\nonumber \\
=&   \frac{1 + \kappa_n}{\kappa_n}  E\left[ \frac{1 + \kappa_n}{\kappa_n} \frac{1}{\gamma} \right] + {\sf Var}(\mu_n)
\nonumber \\
=& \frac{1 + \kappa_n}{\kappa_n}  \frac{\beta_n}{\alpha_n -1 } ,
\label{finalpredvar}
\end{align}
since $\mu$ does not depend on $\gamma$.

\section{Two-Way Bayesian ANOVA}
\label{calcs2wayANOVA}

Here we give the details for working out the three term PPV expansion for
a two-way random effects ANOVA from Subsec. \ref{2wayANOVAeg}. 

\noindent
{\it Step 1:}
Decompose the log-likelihood.  For any $i, j$ we have that
\begin{align}
  ~&\ln p(y_{ij}, \tau_i , \beta_j)\nonumber\\
  =& \ln p(y_{ij} \vert \tau_i, \beta_j) + \ln p(\tau_i) + \ln p(\beta_j)
\nonumber \\
=&
\left[ - \frac{1}{2 \sigma_\epsilon^2} \sum_{i,j} (y_{ij} - \tau_i -\beta_j)^2
- \frac{1}{2 \sigma_\tau^2} \sum_{i} (\tau_i - \tau_0)^2\right.\nonumber\\
&\left.\hspace{.5in}- \frac{1}{2 \sigma_\beta^2} \sum_{j} (\beta_j -\beta_0)^2
\right]+ Extra Terms
\label{LLdecomp}
\end{align}
Apart from the $-1/2$ factor, the part of expression \eqref{LLdecomp}
in square brackets is
\begin{align}
~& \frac{1}{\sigma^2_\epsilon} \sum_i \sum_j \left(y_{ij} +  \tau_i^2 + \beta_j^2\right.\nonumber\\
 &\left.\hspace{.5in} - 2 y_{ij}\tau_i -2 y_{ij}\beta_j -2 \tau_i \beta_j \right)
\nonumber \\
& \quad + \frac{1}{\sigma^2_\tau} \sum_i \left(   \tau_i^2 +\tau_ 0^2   - 2 \tau_i \tau_0  \right)
\nonumber \\
& \quad + \frac{1}{\sigma^2_\beta} \sum_j \left(\beta^2_j + \beta_0^2 - 2 \beta_j \beta_0 \right)
\nonumber \\
=& \sum_j \left( \beta_j^2 \left( \frac{T}{\sigma^2_\epsilon} + \frac{1}{\sigma^2_\beta} \right) 
- 2 \beta_j \left(  \frac{y_{+j} + \tau_+}{\sigma^2_\epsilon} + \frac{\beta_0}{\sigma^2_\beta} \right)
\right)
\nonumber \\
& + \sum_i \left( \tau_i^2 \left( \frac{B}{\sigma^2_\epsilon} + \frac{1}{\sigma^2_\tau} \right)
-
2 \tau_i \left( \frac{y_{i+}}{\sigma^2_\epsilon} + \frac{\tau_0}{\sigma^2_\tau} \right) \right)
\nonumber \\
& + \left(\frac{1}{\sigma^2_\epsilon} \sum_{ij} y_{ij}^2 + \frac{B\beta_0^2}{\sigma^2_\beta}
+ \frac{T\tau_0^2}{\sigma^2_\tau} \right)
\nonumber \\
& \equiv \sum_j (T1)_j + \sum_i(T2)_i + T3 .
\label{terms}
\end{align}

\noindent{\it Step 2:} Use \eqref{terms} to obtain
\begin{eqnarray}
p(\beta_j \vert  {\bf y}, {\bf \tau} ) \sim N\left( \frac{  \frac{y_{+j} + \tau_+} {\sigma^2_\epsilon}  + \frac{\beta_0}{\sigma^2_\beta} }
{  \frac{T}{\sigma^2_\epsilon} + \frac{1}{\sigma^2_\beta} }, 
\left( \frac{T}{\sigma^2_\epsilon} + \frac{1}{\sigma^2_\beta}  \right)^{-1} \right)
\nonumber
\end{eqnarray}
where ${\bf y}$ is the matrix of $y_{ij}$'s, $\tau$ is the vector of $\tau_i$'s, and the subscript $+$ indicates a sum
over the appropriate index.

To see this, set up a completing the square in $\beta_j$.  That is, write
\begin{align}
~&(T1)_j\nonumber\\ 
=&
\beta_j^2 \left( \frac{T}{\sigma^2_\epsilon} + \frac{1}{\sigma^2_\beta} \right) 
-  \frac{ 2 \beta_j \left(  \frac{y_{+j} + \tau_+}{\sigma^2_\epsilon} + \frac{\beta_0}{\sigma^2_\beta} \right)}
{\sqrt{\frac{T}{\sigma^2_\epsilon} + \frac{1}{\sigma^2_\beta}}}
 \sqrt{\frac{T}{\sigma^2_\epsilon} + \frac{1}{\sigma^2_\beta}} 
 \nonumber \\
 & \quad \pm \frac{  \left(  \frac{y_{+j} + \tau_+}{\sigma^2_\epsilon} + \frac{\beta_0}{\sigma^2_\beta} \right)^2  } { \left(  \frac{T}{\sigma^2_\epsilon} + \frac{1}{\sigma^2_\beta}  \right)^2  }
 \label{ET1j}
 \\
 =& 
 \left( \frac{T}{\sigma^2_\epsilon} + \frac{1}{\sigma^2_\beta}  \right)
 \times 
 \left( \beta_j -   \frac{ \frac{y_{+j} + \tau_+}{\sigma^2_\epsilon} + \frac{\beta_0}{\sigma^2_\beta} }
 {  \frac{T}{\sigma^2_\epsilon} + \frac{1}{\sigma^2_\beta}   }  \right)^2 - ET_{1,j}
 \label{squarebeta}
\end{align}
where $ET_{1,j}$ is the positive version of the last term in \eqref{ET1j}.
From \eqref{squarebeta} we get Step 2.

\noindent
{\it Step 3:} Verify that the rest of the `active terms' in the exponent 
\begin{eqnarray}
\sum_i(T2)_i + \sum_j (ET)_{1,j} 
\label{activeterms}
\end{eqnarray}
generate a quadratic form for an appropriate matrix and vector space.
With some foresight, let
\begin{eqnarray}
a = \frac{B}{\sigma^2_\epsilon} + \frac{1}{\sigma^2_\tau} 
\quad 
\hbox{and} 
\quad 
b = -  \frac{1}{ \sigma^4_\epsilon ( \frac{T}{\sigma^2_\epsilon} + \frac{1}{\sigma^2_\beta}) } .
\nonumber
\end{eqnarray}
Also, write
$$
v_i = \left( \frac{y_{i+}}{\sigma^2_\epsilon} + \frac{\tau_0}{\sigma^2_\tau} \right)\left( \frac{  \frac{y_{++}}{\sigma^2_\beta} 
+ \frac{B\beta_0}{\sigma^2_\beta} }{\sigma_\epsilon^2} \right)
$$
and ${\bf v} = (v_1, \ldots, v_T)^T)$.
Now, the active terms from \eqref{activeterms} equal
\begin{align}
 \sum_i &\left( \tau_i^2 \left( \frac{B}{\sigma^2_\epsilon} + \frac{1}{\sigma^2_\tau} \right)
-
2 \tau_i \left( \frac{y_{i+}}{\sigma^2_\epsilon} + \frac{\tau_0}{\sigma^2_\tau} \right) \right)\nonumber\\
-&
\frac{1}{\left( \frac{T}{\sigma^2_\epsilon} + \frac{1}{\sigma^2_\beta}  \right) }
\left[ \sum_j  \left( \frac{y_{+j} + \tau_+}{\sigma^2_\epsilon} + \frac{\beta_0}{\sigma^2_\beta} \right)^2 \right].
\label{expandsq}
\end{align}
The expression in square brackets in \eqref{expandsq} can be re-expressed as
\begin{align}
~&\left[  \frac{B(\sum_i \tau_i^2 + \sum_{k \neq i} \tau_k\tau_i ) }{\sigma_\epsilon^4}\right.\nonumber\\ 
&\left.\hspace{.2in}+ 2 \sum_i \frac{\tau_i(\frac{y_{++}}{\sigma^2_\epsilon} + \frac{B\beta_0}{\sigma^2_\beta} )}{\sigma^2_\epsilon}  
+
\sum_j(OT)_j  \right]
\label{OTj}
\end{align}
where $(OT)_j$ represents the `other terms' in the expansion of the 
expression in square
brackets that do not involve the $\tau_i$'s. For $i=1, \ldots , T$, let
\begin{eqnarray}
w_i = \frac{y_{++}}{\sigma_\epsilon^2} + \frac{\tau_0}{\sigma_\tau^2} 
+
\frac{1}{\sigma_\epsilon^2} \left( \frac{y_{++}}{\sigma_\beta^2} +
\frac{B\beta_0}{\sigma_\beta^2} \right),
\end{eqnarray}
with $W$ meaning the vector of the $w_i$'s.
Using \eqref{OTj} without
$(OT)_j$'s in 
\eqref{expandsq} gives
\begin{align}
\sum_i  & \left[ \tau_i^2 \left(  \left( \frac{B}{\sigma^2_\epsilon} + \frac{1}{\sigma^2_\tau} \right)  
- \frac{ \frac{B}{\sigma^2_\epsilon}}{ \frac{T}{\sigma^2_\epsilon} + \frac{1}{\sigma^2_\beta} } \right) \right.\nonumber\\
&\left. -
2\tau_i \left(
\left( \frac{y_{i+}}{\sigma^2_\epsilon} + \frac{\tau_0}{\sigma^2_\tau} \right)+\left( \frac{  \frac{y_{++}}{\sigma^2_\beta} 
+ \frac{B\beta_0}{\sigma^2_\beta} }{\sigma_\epsilon^2} \right) \right)\right.
\nonumber \\
 &- \left.  \frac{2B} {\sigma^4_\epsilon \left( \frac{T}{\sigma_\epsilon^2} + \frac{1}{\sigma_\beta^2} \right) } 
\tau_i \sum_{k\neq i} \tau_k \right]
\nonumber \\
=& \sum_i \left(  \tau_i^2 (a + Bb) - 2 \tau_i w_i + 2Bb\tau_i \sum_{k \neq i} \tau_k  \right)
\nonumber \\
=&
 \tau^T 
\begin{pmatrix} 
a + Bb & Bb &\ldots &  \ldots & Bb \\
Bb & a+ Bb &\ldots &  \ldots & Bb\\
\vdots & \vdots &\vdots &  \vdots & \vdots \\
Bb & Bb & \ldots & a+Bb & Bb \\
Bb & Bb & \ldots & Bb & a +Bb\\
\end{pmatrix} 
\tau\nonumber\\
&- 2 \tau^T W
\nonumber\\
\equiv& \tau^T A \tau + 2 \tau^T \mu_\tau
\label{innerprod}
\end{align}
where $\mu_\tau  = - W$ and $A$ is of the form
\begin{eqnarray}
A = a I_T + Bb {\bf 1}{\bf 1}^T .
\label{form}
\end{eqnarray}

\noindent{\it Step 4:} Derive the posterior variances and covariances for the $\tau_i$'s.
From the Sherman-Morrison formula we have that
\begin{eqnarray}
A^{-1} = \frac{1}{a} \left(I_T - \frac{Bb  {\bf 1}{\bf 1}^T }{a + Bb T} \right).
\end{eqnarray}
Continuing from \eqref{innerprod} and again completing the square, this part of the exponent in 
the likelihood (see \eqref{LLdecomp}) is
\begin{align}
=& - \frac{1}{2} \left(  \tau^T (A^{-1})^{-1} \tau - 2 \tau \mu_\tau \right)
\nonumber \\
=& - \frac{1}{2} (\tau - \mu_\tau) \sigma^{-1} (\tau - \mu_\tau)\nonumber\\
&+ LowerOrderTerms,\nonumber
\end{align}
for some $n \times n$ matrix $\Sigma$.
Writing $a_{ij}$ for the elements of $A$ and $\sigma_{ij}^{(-1)}$ for the elements in $\Sigma$
we see that for any $i$ and $j$ that
\begin{eqnarray}
a_{ij} \tau_i \tau_j = \sigma^{(-1)}_{ij} \tau_i \tau_j .
\nonumber
\end{eqnarray}
Hence, $A = \Sigma^{-1}$ and $\Sigma = A^{-1}$ and both are symmetric and positive definite.
Now, from the Sherman-Morrison formula we see that
\begin{eqnarray}
{\sf Var}(\tau_i \vert {\bf y}) = \frac{1}{a} \left(1 - \frac{Bb}{a+BbT} \right)
\label{postvar}
\end{eqnarray}
and
\begin{eqnarray}
{\sf Cov}(\tau_i, \tau_j \vert {\bf y}) =  - \frac{Bb}{a(a+BbT)} .
\label{postcovar}
\end{eqnarray}

As a check, we observe that 
\begin{eqnarray}
a + BbT = C \left[ 
\left( \frac{B}{\sigma^2_\epsilon} + \frac{1}{\sigma^2_\tau} \right)  
\left( \frac{T}{\sigma^2_\epsilon} + \frac{1}{\sigma^2_\beta} \right) - \frac{BT}{\sigma^4_\epsilon}  
\right]
\nonumber
\end{eqnarray}
for a suitable $C > 0$ and it is easy to see that the right hand side is strictly positive.  So, \eqref{postvar}
and \eqref{postcovar} are well defined.  Since $b < 0$ both are positive as well.

\noindent
{\it Step 5:}. Now we can derive an expression for the posterior covariance of ${\sf Var}(Y_{ij; n+1} \vert {\bf y})$.
By two uses of the LTV we have
\begin{align}
~&{\sf Var}(Y_{ij; n+1} \vert {\bf y})\nonumber\\ 
=& {\sf E}_{\tau \vert {\bf y}} {\sf E}_{\beta \vert  {\bf y}, \tau } {\sf Var}(Y_{ij; n+1} \vert {\bf y}, \beta, \tau)
\nonumber \\
& + {\sf E}_{\tau \vert {\bf y}}  {\sf Var}_{\beta \vert  {\bf y}, \tau } {\sf E}(Y_{ij; n+1} \vert {\bf y}, \beta, \tau)
\nonumber \\
& + {\sf Var}_{\tau \vert {\bf y}}  {\sf E}_{\beta \vert {\bf y}, \tau } {\sf E}(Y_{ij; n+1} \vert {\bf y}, \beta, \tau) .
\label{2wayANOVAvardecomp}
\end{align}
The first term is
\begin{eqnarray}
{\sf E}_{\tau \vert {\bf y}} {\sf E}_{\beta \vert  {\bf y}, \tau } {\sf Var}(Y_{ij; n+1} \vert {\bf y}, \beta, \tau)
= {\sf E}_{\tau \vert {\bf y}} {\sf E}_{\beta \vert  {\bf y}, \tau }  \left( \sigma^2_\epsilon \right)
= \sigma^2_\epsilon .
\nonumber
\end{eqnarray}
The second term is
\begin{align}
{\sf E}_{\tau \vert {\bf y}}  {\sf Var}_{\beta \vert  {\bf y}, \tau } {\sf E}(Y_{ij; n+1} \vert {\bf y}, \beta, \tau)
=&
{\sf E}_{\tau \vert {\bf y}}  {\sf Var}_{\beta \vert  {\bf y}, \tau } \left( \tau_i + \beta_j \right)\nonumber\\
=&  \frac{1}{ \frac{T}{\sigma_\epsilon^2} + \frac{1}{\sigma^2_\beta} } ,
\end{align}
using the fact that i) $\tau_i$ and $\beta_j$ are independent, ii)
${\sf Var}_{\beta \vert {\bf y}, \tau}(\tau_i) = 0$, and iii)  the result from Step 2.

The third term in \eqref{2wayANOVAvardecomp} is
\begin{align}
~&{\sf Var}_{\tau \vert {\bf y}}  E_{\beta \vert {\bf y}, \tau } E(Y_{ij; n+1} \vert {\bf y}, \beta, \tau) 
\nonumber\\
=&
{\sf Var}_{\tau \vert {\bf y}}  E_{\beta \vert {\bf y}, \tau } \left( \tau_i + \beta_j \right) 
\nonumber \\
=&
{\sf Var}_{\tau \vert {\bf y}}  \left( \tau_i + \frac{ \frac{y_{+j} + \tau_+}{\sigma^2_\epsilon} + \frac{\beta_0}{\sigma^2_\beta}}
{\frac{T}{\sigma^2_\epsilon} + \frac{1}{\sigma^2_\beta} } \right)
\nonumber \\
=& 
\frac
{{\sf Var}_{\tau \vert {\bf y}} \left ( \tau_i \left( \frac{T+1}{\sigma^2_\epsilon} + \frac{1}{\sigma^2_\beta} \right) 
+ \left(\frac{y_{+j}}{\sigma^2_\epsilon} + \frac{\beta_0}{\sigma^2_\beta} \right)
+ \sum_{j \neq i}^T \frac{\tau_j}{\sigma^2_\epsilon}
 \right)}{
\left( \frac{T}{\sigma^2_\epsilon} + \frac{1}{\sigma^2_\beta} \right)^2
} 
 \nonumber \\
=& 
\frac{1}{
\left( \frac{T}{\sigma^2_\epsilon} + \frac{1}{\sigma^2_\beta} \right)^2
} 
 \left\{ 
  \left( \frac{T+1}{\sigma^2_\epsilon} + \frac{1}{\sigma^2_\beta} \right)^2  \cdot \frac{1}{a} \left( 1 - \frac{Bb}{a+bBT}\right)\right.\nonumber\\
  &\left.+
  \frac{(T-1)}{ \sigma^4_\epsilon} \left(\frac{1}{a}\right) \left( 1 - \frac{Bb}{a + bBT} \right)
  \right.
  \nonumber \\
 & +  \left. 
 \frac{2}{\sigma^2_\epsilon} \left( \frac{T+1}{\sigma^2_\epsilon} + \frac{1}{\sigma^2_\beta} \right)  (T-1) 
 \frac{ - Bb}{a(a+ BbT)}
 \right\} .
\label{2wayANOVAexp3}
\end{align}
In the last term in \eqref{2wayANOVAexp3}, we have recognized $2\hbox{Cov}(\tau_i, \tau_j \vert {\bf y})$
and that the number of $\tau_j$'s not equal to a given $\tau_i$ is $T-1$.
\end{appendix}

\begin{acks}[Acknowledgments]
Dustin acknowledges funding from the University of Nebraska
Program of Excellence in Computational Science.  The first author was partially supported by the grant 2413491 from the Division of Mathematical Sciences, National Science Foundation. All authors acknowledge computational support 
from the Holland Computing Center at the University of Nebraska.  
The third author thanks Sujit Ghosh
for helpful discussions about how to use the expansions proposed here.
\end{acks}

\bibliographystyle{imsart-number}
\bibliography{references.bib}

\end{document}

%% file: hmmPic.tex
$
\psmatrix[colsep=.5in,rowsep=.3in,mnode=circle]
Y_1&Y_2&&&&Y_n&Y_{n+1}\\
X_1&X_2&&&&X_n&X_{n+1}\\
\endpsmatrix
\psset{nodesep=3pt,linewidth=1pt,origin={0,0}}
\ncline{->}{1,1}{1,2}
\ncline{->}{1,1}{2,1}
\ncline{->}{1,2}{2,2}
\ncline{->}{1,6}{1,7}
\ncline{->}{1,6}{2,6}
\ncline{->}{1,7}{2,7}
\ncline[linestyle=dotted]{->}{1,2}{1,6}
\ncline[linestyle=dotted]{-}{2,2}{2,6}
\psframe[linestyle=dashed,linearc=14pt,cornersize=absolute,linecolor=blue](-1.5,-1)(.5,1)
\rput(.2,-.55){V_1}
\psframe[linestyle=dashed,linearc=14pt,cornersize=absolute,linecolor=blue](-12,-1)(-2,1)
\rput(-11.5,-.5){V_2}
\psframe[linestyle=dashed,linearc=14pt,cornersize=absolute,linecolor=blue](-11.5,1.25)(-1.75,2.5)
\rput(-2.25,2.25){\mc{D}}
$

%% file: bhm.tex
$
\psmatrix[colsep=.5in,rowsep=.3in]
&&&\circlenode{A}{V_2}&&&\\
&&&\circlenode{A}{V_1}&&&\\
\circlenode{A}{Y_1}&\circlenode{A}{Y_2}&&&&\circlenode{A}{Y_n}&\circlenode{A}{Y_{n+1}}\\
\endpsmatrix
\psset{nodesep=3pt,linewidth=1pt,origin={0,0}}
\ncline{->}{1,4}{2,4}
\ncline{->}{2,4}{3,1}
\ncline{->}{2,4}{3,2}
\ncline{->}{2,4}{3,6}
\ncline{->}{2,4}{3,7}
\ncline[linestyle=dotted]{->}{2,4}{3,3}
\ncline[linestyle=dotted]{-}{3,2}{3,6}
\psframe[linestyle=dashed,linearc=14pt,cornersize=absolute,linecolor=blue](-1.5,-1)(.5,1)
\psframe[linestyle=dashed,linearc=14pt,cornersize=absolute,linecolor=blue](-12,-1)(-2,1)
\rput(-11.5,-.5){\mc{D}}
$

%% file: exmpl4.1.tex
{\Large
$
\psmatrix[colsep=.1\linewidth,rowsep=.05\linewidth]
\mathbf{\pi}&\{1,2,3\}&\{2,1,3\}&\{1,3,2\}&\{3,1,2\}&\{2,3,1\}&\{3,2,1\}\\
\mathbf{T_3}&T^{\{1,2,3\}}_3=0&T^{\{2,1,3\}}_3=0&T^{\{1,3,2\}}_3=0&T^{\{3,1,2\}}_3=0&T^{\{2,3,1\}}_3=0&T^{\{3,2,1\}}_3=0\\[.4\linewidth]
\mathbf{T_2}&T^{\{1,2,3\}}_2=0&T^{\{2,1,3\}}_2=0&T^{\{1,3,2\}}_2=0&T^{\{3,1,2\}}_2=0&T^{\{2,3,1\}}_2=0&T^{\{3,2,1\}}_2=0\\[.2\linewidth]
\mathbf{T_1}&T^{\{1,2,3\}}_1=0&&T^{\{1,3,2\}}_1=0&T^{\{3,1,2\}}_1=0&&T^{\{3,2,1\}}_1=0\\[.075\linewidth]
\mathbf{T_1}&&T^{\{2,1,3\}}_1=0&&&T^{\{2,3,1\}}_1=0&
\endpsmatrix
\psset{nodesep=3pt,linewidth=.5pt,origin={0,0},doubleline=true}
\ncline{<->}{2,2}{2,3}
\ncline{<->}{2,4}{2,5}
\ncline{<->}{2,6}{2,7}
\ncline{->}{2,2}{3,4}
\ncline{->}{2,2}{3,6}
\ncline{->}{2,3}{3,4}
\ncline{->}{2,3}{3,6}
\ncline[linestyle=dashed]{->}{2,4}{3,2}
\ncline[linestyle=dashed]{->}{2,4}{3,7}
\ncline[linestyle=dashed]{->}{2,5}{3,2}
\ncline[linestyle=dashed]{->}{2,5}{3,7}
\ncline[linestyle=dotted]{->}{2,6}{3,3}
\ncline[linestyle=dotted]{->}{2,6}{3,5}
\ncline[linestyle=dotted]{->}{2,7}{3,3}
\ncline[linestyle=dotted]{->}{2,7}{3,5}
\ncline{<->}{4,2}{4,4}
\ncline{<->}{4,5}{4,7}
\ncline{<->}{5,3}{5,6}
\ncline{->}{3,4}{4,5}
\ncline{->}{3,4}{4,7}
\ncline{->}{3,6}{4,5}
\ncline{->}{3,6}{4,7}
\ncline[linestyle=dashed]{->}{3,2}{5,3}
\ncline[linestyle=dashed]{->}{3,7}{5,3}
\ncline[linestyle=dashed]{->}{3,2}{5,6}
\ncline[linestyle=dashed]{->}{3,7}{5,6}
\ncline[linestyle=dotted]{->}{3,3}{4,2}
\ncline[linestyle=dotted]{->}{3,5}{4,2}
\ncline[linestyle=dotted]{->}{3,3}{4,4}
\ncline[linestyle=dotted]{->}{3,5}{4,4}
$
}